\documentclass[preprint,final,3p,times,12.8pt,authoryear]{elsarticle}
\usepackage{mathrsfs}
\usepackage{amsthm,amsmath,amssymb}  
\usepackage{algorithm,algorithmicx,algpseudocode} 
\usepackage{array}
\usepackage{multirow}
\usepackage{bm}
\usepackage{hyperref}
\hypersetup{hypertex=true,
            colorlinks=true,
            linkcolor=blue,
            anchorcolor=blue,
            citecolor=blue}

\newtheorem{theorem}{Theorem}[section]
\newtheorem{lemma}[theorem]{Lemma}
\newtheorem{assumption}[theorem]{Assumption}

\newtheorem{remark}{Remark}[section]

\numberwithin{equation}{section}
\numberwithin{figure}{section}

\makeatletter
\newenvironment{breakablealgorithm}
  {
   \begin{center}
     \refstepcounter{algorithm}
     \hrule height.8pt depth0pt \kern2pt
     \renewcommand{\caption}[2][\relax]{
       {\raggedright\textbf{\ALG@name~\thealgorithm} ##2\par}%
       \ifx\relax##1\relax 
         \addcontentsline{loa}{algorithm}{\protect\numberline{\thealgorithm}##2}%
       \else 
         \addcontentsline{loa}{algorithm}{\protect\numberline{\thealgorithm}##1}%
       \fi
       \kern2pt\hrule\kern2pt
     }
  }{
     \kern2pt\hrule\relax
   \end{center}
  }
\makeatother 

\journal{Computers and Mathematics with Applications}

\begin{document}

\begin{frontmatter}

\title{Low rank approximation method for perturbed linear systems with applications to elliptic type stochastic PDEs}

\author{Yujun Zhu\fnref{addr0}}
\ead{yujun_zhu@hust.edu.cn}

\author{Ju Ming\fnref{addr0}\corref{cor1}} 
\ead{jming@hust.edu.cn}
\cortext[cor1]{Corresponding author}

\author{Jie Zhu\fnref{addr0}}
\ead{zj2021@hust.edu.cn}

\author{Zhongming Wang\fnref{addr1}}
\ead{zwang6@fiu.edu}

\address[addr0]{School of Mathematics and Statistics, Huazhong University of Science and Technology, Wuhan, China}

\address[addr1]{Department of Mathematics and Statistics, Florida International University, Miami, FL, USA}

\begin{abstract}
{In this paper, we propose a low rank approximation method for efficiently solving stochastic partial differential equations. Specifically, our method utilizes a novel low rank approximation of the  stiffness matrices, which can significantly reduce the computational load and storage requirements associated with matrix inversion without losing accuracy. To demonstrate the versatility and applicability of our method, we apply it to address two crucial uncertainty quantification problems: stochastic elliptic  equations and  optimal control problems governed by stochastic elliptic PDE constraints. Based on varying dimension reduction ratios, our algorithm exhibits the capability to yield a high precision numerical solution for stochastic partial differential equations, or provides a rough representation of the exact solutions as a pre-processing phase. Meanwhile, our algorithm for solving stochastic optimal control problems allows a diverse range of gradient-based unconstrained optimization methods, rendering it particularly appealing for computationally intensive large-scale problems. Numerical experiments are conducted and the results provide strong validation of the feasibility and effectiveness of our algorithm.
}
\end{abstract}


\end{frontmatter}

\section{Introduction}

The systems of linear equations have been widely used in many real-life fields such as natural science (\citealp{jin2020saturation}; \citealp{barz2013solving}), engineering (\citealp{bellman2016introduction}; \citealp{dehghan2012fourth}), economics (\citealp{poongodi2020prediction}), and industry (\citealp{jin2018noise}; \citealp{ringot2007vitro}; \citealp{sonneveld2009idr}). In most of these applications, the coefficient matrices in the linear equations are determined by physical measurements or approximations. Such measurement or approximation errors  could be considered as random perturbations to the coefficient matrices \citep{chen2002perturb}. Due to the uncertainties in the coefficient matrices, it is no longer feasible to obtain the inverse of the perturbed coefficient matrices. Therefore we will focus on the statistical information of the system output, i.e., the statistical moments.  Adopting  the Monte Carlo method \citep{metropolis1949monte} for discretizing stochastic space, one can transform the problem to solve $M$ sampling linear equations with similar coefficient matrices. Considering a collection of  $M\times M$ coefficient matrices $\{\mathbb{A}_m\}^{M}_{m=1}$ and the right hand vector $\bm{b}$, our goal is to solve the following linear equations:

\begin{equation} \label{eq:pls}
\mathbb{A}_m \bm{x}_m = \bm{b}, \enspace m = 1,...,M.
\end{equation}

However, we will encounter expensive computational  and memory  costs for $\{\mathbb{A}_m\}^{M}_{m=1}$, when the amount of data required to obtain reliable statistical moments and degrees of freedom of the coefficient matrices are large. To reduce both costs,  it is natural to consider the dimension reduction techniques for the matrices. The aim of dimensional reduction is to obtain a lower dimensional compact data representation with little loss of information. One of the most commonly used dimensional reduction model is vector space model \citep{turk1991eigenfaces, zhao2003face}, which  is based on vectors selection. Within this model, each coefficient matrix is treated as a collection of column vectors. Many algorithms have been proposed based on this model in various applications, such as face recognition \citep{turk1991eigenfaces}, machine learning \citep{castelli2003csvd} and information retrieval \citep{berry1995using}. A well-known technique based on this vector space model  is the low rank approximation by using singular value decomposition (SVD). An appealing property of low rank approximation of matrices (LRAM) via SVD is that it can achieve the smallest reconstruction error among all approximations with the constrain of same rank in Euclidean distance \citep{eckart1936approximation}. 

However, the traditional SVD-based LRAM confronts practical computational limits due to the large time and space complexities in dealing with large matrices. Many attempts are made to conquer the computational cost obstacle. For instance, a generalized LRAM method, proposed in \citep{ye2004generalized}, is proved to have less computational time than the traditional SVD-based method in practical applications. To further lessen computation cost, a simplified generalized LRAM method is proposed in \citep{lu2008simplified}, which further simplifies the structure and the projection matrices.  The robust generalized LRAM algorithm does well in the presence of large sparse noise or outliers \citep{shi2015robust, zhao2016robust}.  The non-iterative LRAM method in \citep{liu2006non} uses an analytical rather form than an iterative manner.  Overall, the current research mainly focuses on the data compression for large-scale matrices, and there are few literature exploring the application of LRAM in approximating the inverses of matrices or solving the perturbed linear systems.

There is a special case of random linear systems with perturbations in the coefficient matrices: the probabilistic discretized formulation of stochastic partial differential equations (SPDEs), which depict many physical and engineering models involving uncertain data or parameters. Due to the wide application in the areas of applied science and engineering, both academia and industry have shown growing interest in designing efficient numerical methods for solving SPDEs \citep{babuvska2002solving, babuvska2003solving, allen1998finite, burkardt2007reduced, du2002numerical, gunzburger1996finite}. 

\begin{figure}[ht]
\centering
\includegraphics[width=0.9\textwidth]{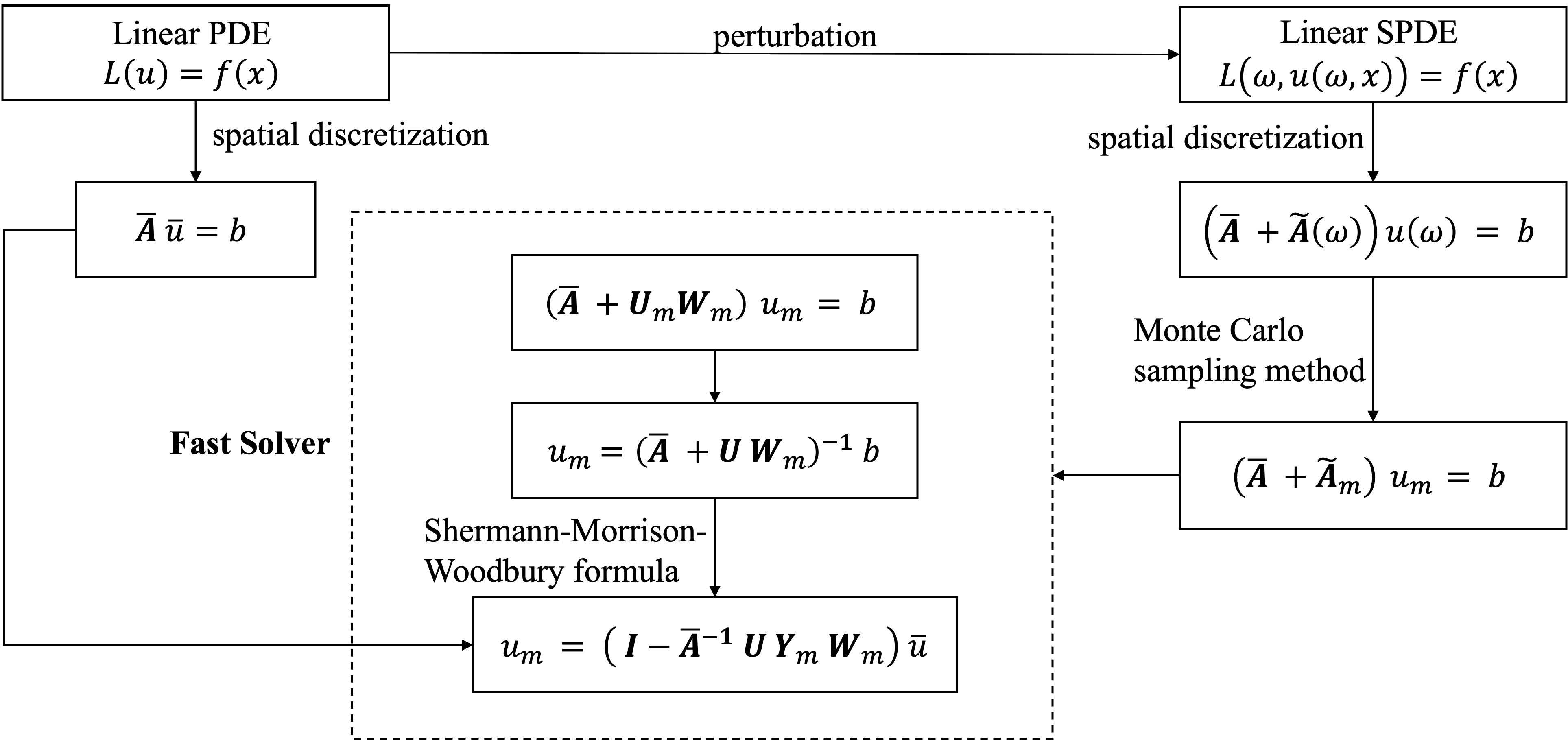}
\caption{The flowchart of the fast solver algorithm, where matrices $\bm{U}, \bm{W}_m$ are the low rank approximation of the perturbed stiffness matrices $\widetilde{\mathbb{A}}_m$.}
\label{fig1.1}
\end{figure}

In this article, we consider the linear partial differential equations with perturbed inputs $\omega$, i.e.,

$$\mathcal{L}(\omega, \bm{u}(\bm{x},\omega)) = f(\bm{x}),$$

where $L$ is a linear differential operator and $f(\bm{x})$ is a smooth function. The corresponding spatial discretization is 

$$\mathbb{A}(\omega) \bm{u}(\omega) = \bm{b},$$ 

where uncertainties lie in the stiffness matrices $\mathbb{A}(\omega)$.  The stiffness matrices $\mathbb{A}(\omega)$ are composed by two parts, the initial deterministic matrices $\overline{\mathbb{A}}$ and the stochastic matrices $\widetilde{\mathbb{A}}(\omega)$ stemmed from perturbations.
We then discretize the probabilistic space by the Monte Carlo (MC) sampling method and obtain the linear system in the form \eqref{eq:pls} 

$$\mathbb{A}_m \bm{u}_m = \bm{b},$$ 

where $\mathbb{A}_m = \overline{\mathbb{A}} + \widetilde{\mathbb{A}}_m$. Note that $\overline{\mathbb{A}}$ is deterministic and invariant to $m$.  We then propose a novel low rank approximation method for the collection of large-scale matrices $\{\widetilde{\mathbb{A}}_m\}_{m=1}^M$, which is derived from an SPDE rather than general noise. Finally we could obtain the numerical solutions by the Shermann-Morrison-Woodbury formula. Figure \ref{fig1.1} presents the schematic flowchart of our algorithm. 

Our algorithm can significantly reduce the computational complexity and storage requirement. In addition, the numerical solutions rely on the initial solution $\bm{x}_0 = \mathbb{A}^{-1} \bm{b}$, and therefore we could make use of $\bm{x}_0$ and obtain the new solution $\bm{x}_m$ directly rather than to solve a brand new linear system. Our algorithm is applied  to stochastic elliptic partial differential equation and stochastic optimal control problems with elliptic PDE constraint. Depending on different dimension reduction ratio $\tau$, our algorithm can construct a high precision numerical solution or roughly depict the main sketch of the exact solution as a pre-processing. Moreover, the algorithm can make good use of various types of gradient-based unconstrained optimization methods, which makes it more attractive for large-scale problems from a computational point-of-view. Numerical results also validate the feasibility and the effectiveness of the proposed algorithms. 

The remainder of this article is organized as follows. Section 2 presents a brief overview of some related function spaces and notations. In Section 3, we introduce the existing low rank approximation methods and propose a novel low rank matrix approximation approach with less reconstruction error and CPU memory. Then we propose a fast algorithm for efficiently solving the perturbed linear equations and analyze its computational complexity and storage. The applications to the stochastic elliptic partial differential equation and the stochastic optimal control problem with the elliptic SPDE constrain are studied in Section 4. This section also includes the error analysis of the fast numerical methods and their numerical performance studies. Conclusions and discussions of future work can be found in Section 5.

\section{Preliminaries}

We begin by recalling some required function spaces and notations. Throughout this article, we use the standard notations for Sobolev spaces \citep{adams2003sobolev}. Let $L^p(D),  1 \leq p \leq \infty$, denote usual Lebesgue space on smooth domain $D\in \mathbb{R}^n$; $\Vert \cdot \Vert = \Vert \cdot \Vert_{L^2(D)}$ denote the $L^2(D)$-norm induced by the inner product $\langle f,g \rangle = \int_D fg dx, \, \forall f,g\in L^2(D)$. $H^r(D)$ is a Sobolev space for all real numbers $r$ with norms $\Vert y \Vert_{H^r(D)}$ and semi-norms $\vert y \vert_{H^r(D)}$, where

\begin{equation} \label{eq2.1}
\Vert y \Vert_{H^r(D)} = \sum_{\vert \bm{\alpha} \vert \leq r} \Vert \frac{\partial^{\bm{\alpha}} u}{\partial u^{\bm{\alpha}}} \Vert^2_{L^2(D)},
\end{equation}

\noindent and

\begin{equation} \label{eq2.2}
\vert y \vert_{H^r(D)} = \sum_{\vert \bm{\alpha} \vert = r} \Vert \frac{\partial^{\bm{\alpha}} u}{\partial u^{\bm{\alpha}}} \Vert^2_{L^2(D)}.
\end{equation}
\noindent Here $\alpha$ is a multi-index with non-negative integer components $\{\alpha_i\}$, and $\bm{\alpha} = \sum_i \alpha_i$.

Sobolev spaces have forms

\begin{equation} \label{eq2.3}
H^1(D) = \{ y \in L^2(D), \partial_{x_i} y \in L^2(D), i = 1,...,n\},
\end{equation}

\noindent and

\begin{equation} \label{eq2.4}
H_0^1(D) = \{ y \in H^1(D), y \vert_{\partial D} = 0\}.
\end{equation}

\noindent Clearly, $H_0^1(D)$ is a subspace of $H^1(D)$. Moreover, we define $\Vert y \Vert_{H_0^1(D)} = \Vert \nabla y \Vert$, and $H^{-1}(D)$  is the dual space for $H_0^1(D)$.

Let $(\Omega, \mathscr{F}, \mathbb{P})$ be a complete probability space, where $\Omega$ denotes the set of outcomes, $\mathscr{F} \subset 2^\Omega$ is the $\sigma$-algebra of events, and $\mathbb{P}: \mathscr{F} \rightarrow [0,1]$ is a complete probability measure, respectively. If $X$ is a real random variable in $(\Omega, \mathscr{F}, \mathbb{P})$, then its expectation is given by

\begin{equation} \label{eq2.5}
\mathbb{E}[X] = \int_\Omega X(\omega) \mathbb{P}(d \omega) = \int_{\mathbb{R}^n} x \rho(x) dx,
\end{equation}

\noindent where $\rho$ is a joint PDF for $X$, defined on a Borel set $\mathcal{B}$ of $\mathbb{R}$, that is, $\rho(\mathcal{B}) = \mathbb{P}(X^{-1}(\mathcal{B}))$.

Define the stochastic Sobolev space as

\begin{equation} \label{eq2.6}
L^2(\Omega; H^r(D)) = \{ y: D \times \Omega \rightarrow \mathbb{R} \enspace \vert \enspace \Vert y \Vert_{L^2(\Omega; H^r(D))} \, < \, \infty \},
\end{equation}

\noindent equipped with the norm 

\begin{equation} \label{eq2.7}
\Vert y \Vert_{L^2(\Omega; H^r(D))} = \int_\Omega \Vert y \Vert_{H^r(D)} d \mathbb{P} = \mathbb{E}[\Vert y \Vert_{H^r(D)}].
\end{equation}

The stochastic Sobolev space $L^2(\Omega; H^r(D))$ is a Hilbert space, and it is isomorphic to the tensor product space $L^2(\Omega) \otimes L^2(H^r(D))$. For simplicity, we set 

\begin{equation} \label{eq2.8}
\mathcal{H}^r(D) = L^2(\Omega; H^r(D)), \mathcal{H}_0^r(D) = L^2(\Omega; H_0^r(D)), 
\mathcal{L}^2(D) = L^2(\Omega; L^2(D)).
\end{equation}

\section{Low Rank Approximation Method for Perturbed Linear Systems}

\subsection{Problem setting and perturbation splitting}

Consider a PDE operator $\mathcal{L}$ and suppose the force $f(x): D \rightarrow \mathbb{R}$ is a deterministic term, then a typical stochastic partial differential equation can be stated as follows: find $u(x, \omega): \Omega \times D \rightarrow \mathbb{R}$ satisfying

\begin{equation} \label{eq3.1}
\mathcal{L}(\omega, u(x, \omega)) = f(x), \enspace a.e. \enspace in \enspace D, 
\end{equation}

\noindent over a bounded, Lipschitz domain $D \subset \mathbb{R}^d$, $d = 1,2,3$ and a complete probability space $(\Omega, \mathscr{F}, \mathbb{P})$, and equipped with suitable boundary conditions. The finite dimensional random variable $\omega$ has a joint probability density function (PDF) $\rho: \Omega \rightarrow \mathbb{R}^+$, with $\rho \in L^\infty(\Omega)$. In order to predict statistical behaviors of the physical system in Eq. (\ref{eq3.1}), our goal is to obtain the approximation of the multi-dimensional statistical quantities of interest (QoI) \citep{smith2013uncertainty}

\begin{equation} \label{eq3.2}
\mathbb{E}[u](x) = \int_\Omega u(x, \omega) \rho(\omega) d\omega, \enspace where \enspace \omega \in \Omega \enspace and \enspace x \in D. 
\end{equation}

It is generally infeasible to obtain analytic solutions of the SPDE in Eq. (\ref{eq3.1}), and thereby we discretize the equation spatially and approximate it by a random matrix problem as read 

\begin{equation} \label{eq3.3}
\mathbb{A}(\omega) \: \bm{u}(\omega) \: = \: \bm{b}.
\end{equation}

To tackle with the uncertainties in the inherited stiff  coefficient matrix $\mathbb{A}(\omega)$, the probabilistic space also needs to be approximated by a finite-dimensional space. The stochastic matrices $\mathbb{A}(\omega)$ arising from spatial discretization are generally large and sparse, which results in a high-dimensional problem.  The Monte Carlo (MC) method is a natural choice for numerical implementation to lessen the curse of dimensionality, 

$$\mathbb{A}_m \bm{u}_m = \bm{b},\quad m=1,\cdots,M.$$ 
 
However, the precision of the MC method is achieved only by a sufficient large amounts of samples $M$, which again leads to large amounts of matrices inversions. Direct inversion of all  MC realization matrices $\mathbb{A}_m$ requires both huge computational complexity and memory storage.  Observing the similarities between the matrices $\mathbb{A}_m$ of all MC samples, we propose to perturb  $A_m$ into  the form of 

\begin{equation}\label{eq3.4}
\mathbb{A}_m=\overline{\mathbb{A}} \: + \: \widetilde{\mathbb{A}}(\omega),
\end{equation}

where $\overline{\mathbb{A}}$ is the initial deterministic stiffness matrix, and the stochastic stiffness matrix $\widetilde{\mathbb{A}}(\omega)$ stemmed from perturbations, which generally has a relatively low rank. Then we apply the MC sampling method to discretize the probability space $\Omega$ and obtain the following system of equations.

\begin{equation} \label{eq3.5}
(\overline{\mathbb{A}} \: + \: \widetilde{\mathbb{A}}_m) \: \bm{u}_m \: = \: \bm{b}, \enspace m = 1,...,M,
\end{equation}

where $\{\widetilde{\mathbb{A}}_m\}_{m=1}^M$ denote the MC realizations. Note that $\overline{\mathbb{A}}$ is fixed and MC sampling is only needed for the low-ranked $\widetilde{\mathbb{A}}(\omega)$.

The QoI in Eq. (\ref{eq3.2}) is now approximated by

\begin{equation} \label{eq3.6}
\bm{u}_{approx} \: \approx \: \mathbb{E}[u](x) \: = \: \frac{1}{M} \sum_{m=1}^M \bm{u}_m.  
\end{equation}

\begin{remark}
In Figure \ref{fig3.1}, we show the stiffness matrices of  $\overline{\mathbb{A}}$ and $\mathbb{A}_m, m=458, 154$ for the perturbed system of discretized elliptic PDE with diffusion coefficients. We do observe the similarity between the stiffness matrices and thus it is feasible to explore a generalized approximate formulation for the collection of matrices $\{\widetilde{\mathbb{A}}_m\}_{m=1}^M$, which will reduce the storage cost. Meanwhile, the generalized matrix formulation could be also applied in solving the PDEs with perturbations. We do not need to deal with the perturbed systems by directly solving each MC realization of the corresponding linear equations and thereby decrease the computational complexity.

\begin{figure}[ht]
\centering
\includegraphics[width=1\textwidth]{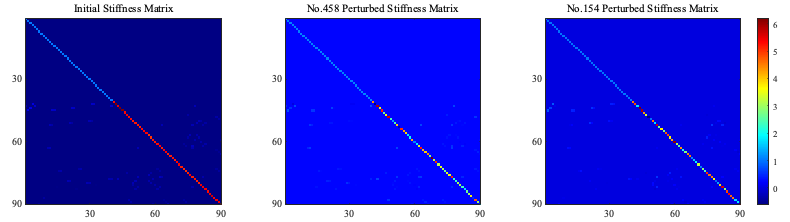}
\caption{Visualization of the stiffness matrices of the elliptic PDE with random diffusion coefficients.}
\label{fig3.1}
\end{figure}

\end{remark}

\begin{remark}
The matrix splitting \eqref{eq3.4} could also take advantage of the initial numerical solution $\overline{\bm{u}} = \overline{\mathbb{A}}^{-1} \bm{f}$ and the low-rank structure of $\widetilde{\mathbb{A}}(\omega)$. Meanwhile, from the aspect of signal processing and deep learning \citep{hemanth2017deep}, the pre-processing of subtracting the statistical mean, known as data normalization, removes the common parts in the MC samples and highlight the individual differences, and thus improving the generalization ability of our algorithm. 
\end{remark}

\subsection{Low Rank Approximation}

In order to obtain an accurate numerical solution in Eq. (\ref{eq3.6}), we need to solve a large-scale linear system in Eq. (\ref{eq3.5}) with a collection of  matrices $\{\widetilde{\mathbb{A}}_m\}_{m=1}^M$.  Since $\overline{\mathbb{A}}$ is independent of sample size $M$,  it is important to explore a high-performance algorithm for dealing with the stochastic 

$$ \widetilde{\mathbb{A}}_m \: \bm{u}_m \: = \: \bm{b}, \enspace m = 1,...,M.$$

To lessen the expensive storage requirement for the perturbed stiffness matrices $\{\widetilde{\mathbb{A}}_m\}_{m=1}^M$, we turn to the low rank approximation technique, which aims to reduce the dimensionality and obtain more compact representations of the original matrices with limited loss of information.

In this section, we develop a SVD-based low rank approximation method to deal with the large-scale matrices $\{\widetilde{\mathbb{A}}_m\}_{m=1}^M$. The goal of low rank approximation is to approximate a large-scale matrix $\mathbb{B}$ with a lower-rank alternative $\widetilde{\mathbb{B}}$, so that solving a given problem on $\widetilde{\mathbb{B}}$ gives a good approximation to the solution on $\mathbb{B}$. A lower-rank sketch means fewer degrees of freedom and less redundancy, and thereby significantly reducing the run-time and memory storage. Mathematically, the optimal rank-$k$ approximation of a matrix $\mathbb{B} \in \mathbb{R}^{N \times N}$ under the Frobenius norm is formulated as a rank-constrained minimization problem (\ref{eq3.7}):  find a matrix $\widetilde{\mathbb{B}} \in \mathbb{R}^{N \times k}$ such that

\begin{equation} \label{eq3.7}
\widetilde{\mathbb{B}}^* \: = \: \mathop{\arg\min}\limits_{rank(\widetilde{\mathbb{B}}) \; = \; k} \enspace \Vert \mathbb{B} - \widetilde{\mathbb{B}} \Vert_F. 
\end{equation}

\subsubsection{Traditional SVD}
The optimal rank-$k$ approximation $\widetilde{\mathbb{B}}^*$ admits an analytical solution in terms of the truncated singular value decomposition (TSVD) of the matrix $\mathbb{B}$, as stated in the following theorem. 

\begin{theorem}[\citealp{eckart1936approximation}] \label{th3.1} 

Let $\mathbb{B} = \mathbb{U} \Sigma \mathbb{V}^T \in \mathbb{R}^{N \times N}$ be the SVD of $\mathbb{B}$, and let $\mathbb{U}, \Sigma$ and $\mathbb{V}$ partitioned as follows:

\begin{equation} 
\mathbb{U} \: := \: \begin{bmatrix} \mathbb{U}_1 & \mathbb{U}_2 \end{bmatrix}, \quad \Sigma \: := \: \begin{bmatrix} \Sigma_1 & 0 \\ 0 & \Sigma_2 \\ \end{bmatrix}, \quad and \quad \mathbb{V} \: := \: \begin{bmatrix} \mathbb{V}_1 & \mathbb{V}_2 
\nonumber
\end{bmatrix},
\end{equation}

\noindent where $\mathbb{U}_1, \mathbb{V}_1 \in \mathbb{R}^{N \times k}$ and $\Sigma_1 \in \mathbb{R}^{k \times k}$. Then the rank-k matrix, obtained from the TSVD, $\widetilde{\mathbb{B}^*} = \mathbb{U}_1 \Sigma_1 \mathbb{V}_1^T$, satisfies that

\begin{equation} \label{eq3.8}
\Vert \mathbb{B} - \widetilde{\mathbb{B}}^* \Vert_F \: = \: \mathop{\min}\limits_{rank(\widetilde{\mathbb{B}}) \leq k} \: \Vert \mathbb{B} - \widetilde{\mathbb{B}} \Vert_F.
\end{equation}

\noindent The minimizer $\widetilde{\mathbb{B}}^*$ is unique if and only if $\sigma_{k+1} \neq \sigma_k$.

\end{theorem}

Theorem \ref{th3.1} reveals the relationship between the low rank approximation and the SVD technique, which guarantees that such approximation is optimal in terms of minimizing the Frobenius norm of the difference among all approximations with the constraint of same rank in Euclidean distance. 

\subsubsection{General low rank approximation of matrix (GLRAM)}
However, due to the expensive SVD computation, the application of the technique in large matrices encounters practical limits both in time and space aspect. Therefore, the generalized low rank approximations of matrices method \citep{ye2004generalized} is proposed to alleviate the high SVD computational cost, which aims to find two matrices $\mathbb{L}, \mathbb{R} \in \mathbb{R}^{N \times k}$ with orthonormal columns, such that 

\begin{equation} \label{eq3.9}
\mathop{\min}\limits_{\mathbb{L}^T \mathbb{L} = I_k , \mathbb{R}^T \mathbb{R} = I_k} \quad \sum_{m=1}^M \Vert \widetilde{\mathbb{A}}_m - \mathbb{L} \mathbb{M}_m \mathbb{R}^T \Vert_F^2. 
\end{equation}

\begin{algorithm}[h]
\caption{The generalized low rank approximations of matrices algorithm} \label{alg1}
\begin{algorithmic}[1]
	
\Require
Matrices $\{\widetilde{\mathbb{A}}_m\}_{m=1}^M$, and the dimension $k$
\Ensure
Matrices $\mathbb{L}, \mathbb{R}$ and $\{\mathbb{M}_m\}_{m=1}^M$

\State Obtain initial $\mathbb{L}_0$ and set $i = 1$.

\While{not convergent}

\State Form matrix $\mathbb{N}_R(\mathbb{L}) = \sum_{m=1}^M \widetilde{\mathbb{A}}_m^T \mathbb{L}_{i-1} \mathbb{L}_{i-1}^T \widetilde{\mathbb{A}}_m$.

\State Compute the $k$ eigenvectors $\{\phi_R^j\}_{j=1}^k$ of $\mathbb{N}_R$ corresponding to the largest $k$ eigenvalues, and let $\mathbb{R}_{i} = \left[ \phi_R^1,...,\phi_R^k \right]$.

\State Form matrix $\mathbb{N}_L(\mathbb{R}) = \sum_{m=1}^M \widetilde{\mathbb{A}}_m \mathbb{R}_{i} \mathbb{R}_{i}^T \widetilde{\mathbb{A}}_m^T$.

\State Compute the $k$ eigenvectors $\{\phi_L^j\}_{j=1}^k$ of $\mathbb{N}_L$ corresponding to the largest $k$ eigenvalues, and let $\mathbb{L}_{i} = \left[ \phi_L^1,...,\phi_L^k \right]$.

\EndWhile

\State \Return Matrices $\mathbb{L} = \mathbb{L}_{i}, \mathbb{R} = \mathbb{R}_{i}, \mathbb{M}_m = \mathbb{L}^T \widetilde{\mathbb{A}}_m \mathbb{R}$.

\end{algorithmic}
\end{algorithm}

The stopping criteria for Line 2 - 8 is given by

\begin{equation} \label{eq3.10}
RMSRE(M) \: := \: \sqrt{\frac{1}{M} \sum_{m=1}^M \Vert \widetilde{\mathbb{A}}_m - \mathbb{L} \mathbb{M}_m \mathbb{R}^T \Vert_F^2},
\end{equation}

\noindent where RMSRE stands for the root mean square reconstruction error. However, it is noteworthy that the convergence of this iterative algorithm cannot be guaranteed theoretically. Meanwhile, the GLRAM method is mainly applied in computer vision and signal processing, which have different requirement for the matrix reconstruction errors from the numerical computation in PDE settings. Therefore, we expect to improve such dimensionality reduction technique such that both the data compression ratio, computational complexity and approximation accuracy for matrices are taken into account, and apply the technique to solve the perturbed linear system in Eq. (\ref{eq3.5}).

\subsection{Proposed LRAM-based Methods}

We begin by obtaining the optimal rank-$k$ approximation of the collection of the perturbed stiffness matrices. Compared to the traditional SVD and GLRAM method, we propose to apply a different data representation of matrices which allows better dimension reduction and data compression. The optimal rank-$k$ approximation of matrices $\{\widetilde{\mathbb{A}}_m\}_{m=1}^M$ is stated as below: find optimal matrices $\mathbb{U}, \{\mathbb{V}_m\}_{m=1}^M \in \mathbb{R}^{N \times k}$ such that

\begin{equation} \label{eq3.11}
\mathop{\min}\limits_{\mathbb{U} \in \mathbb{R}^{N \times k}: \mathbb{U}^T \mathbb{U} = I_k \atop \mathbb{V}_m \in \mathbb{R}^{N \times k}: m = 1, 2,...,M} \quad  \sum_{m=1}^M \Vert \widetilde{\mathbb{A}}_m - \mathbb{U} \mathbb{V}_m^T \Vert_F^2, 
\end{equation}

We now present our algorithm to solve the minimization problem in Eq. (\ref{eq3.11}). We first introduce the following lemma needed for the our main result.

\begin{lemma}\label{le3.2}
Let $\mathbb{U} \in \mathbb{R}^{N \times k}$ and $\mathbb{N} \in \mathbb{R}^{N \times N}$ be symmetric, then the following optimization problem:

\begin{equation}
\begin{split}
\max & \quad tr(\mathbb{U}^T \mathbb{N} \mathbb{U}),\\
s.t. & \quad \mathbb{U}^T \mathbb{U} = I_k, 
\nonumber
\end{split}
\end{equation}

\noindent has the optimal solution satisfying that the columns of $\mathbb{U}$ are eventually the first $k$th eigenvectors of $\mathbb{N}$.

\end{lemma}

\begin{proof}

Let $\mathbb{Q} = [\mathbb{U} \: \mathbb{P}]$ be an orthogonal matrix, i.e. $\mathbb{U} = \mathbb{Q} \begin{bmatrix} \mathbb{I}_k \\ 0 \end{bmatrix}$. Then we have

\begin{equation}
\begin{split}
tr(\mathbb{U}^T \mathbb{N} \mathbb{U}) \: &= \: tr(\; [\mathbb{I}_k 0] \; \mathbb{Q}^T \mathbb{N} \mathbb{Q} \; \begin{bmatrix} \mathbb{I}_k \\ 0 \end{bmatrix}) \: = \: tr(\; \mathbb{Q}^T \mathbb{N} \mathbb{Q} \; [\mathbb{I}_k 0] \; \begin{bmatrix} \mathbb{I}_k \\ 0 \end{bmatrix})\\
&= \: tr(\; \mathbb{Q}^T \mathbb{N} \mathbb{Q} \; [\mathbb{I}_k 0] \; \begin{bmatrix} \mathbb{I}_k & 0 \\ 0 & 0 \end{bmatrix}) \: = \: \sum_{i=1}^k (\mathbb{Q}^T \mathbb{N} \mathbb{Q})_{i,i} =: \: \sum_{i=1}^k \mathbb{G}_{i,i}.
\nonumber 
\end{split}
\end{equation}

\noindent Thus the problem is turned to find an orthogonal matrix $\mathbb{Q}$ to maximize $\sum_{i=1}^k (\mathbb{Q}^T \mathbb{N} \mathbb{Q})_{i,i}$. Let $\mathbb{N} = \mathbb{E} \Sigma \mathbb{E}^T$ be the eigen-decomposition of $\mathbb{N}$, where $\Sigma = diag(\lambda_1, ..., \lambda_N)$, $\lambda_1 \geq ... \geq \lambda_N$ are the eigenvalues of $\mathbb{N}$ in descending order and $\mathbb{E}$ is an orthogonal matrix consisting of eigenvectors of $\mathbb{N}$. We have $\mathbb{G} = \mathbb{Q}^T \mathbb{N} \mathbb{Q} = \mathbb{Q}^T \mathbb{E} \Sigma \mathbb{E}^T \mathbb{Q}$, and $\lambda_1, ..., \lambda_N$ are also eigenvalues of $\mathbb{G}$. Let $d := diag(\mathbb{G})$ and $\lambda := diag(\Sigma)$, then by the Schur-Horn theorem \citep{horn1954doubly}, $d$ is majored by $\lambda$. As a result, $\sum_{i=1}^k \mathbb{G}_{i,i} \leq \sum_{i=1}^k \Sigma_{i,i}$ and $\mathbb{E}^T \mathbb{U} = \mathbb{I}_k$, i.e. $\mathbb{U} = \mathbb{E}$ is consisting of the first $k$ eigenvectors of $\mathbb{N}$.

\end{proof}

We are now ready to state our main results in Theorem \ref{th3.3} and Theorem \ref{th3.4} for finding optimal low-rank matrices $\mathbb{U}$ and $\mathbb{V}_m$.

\begin{theorem}\label{th3.3}

Let $\mathbb{U}$ and $\{\mathbb{V}_m\}_{m=1}^M$ be the optimal solution to the minimization problem in Eq. (\ref{eq3.11}), then $\mathbb{V}_m = \widetilde{\mathbb{A}}_m^T \mathbb{U}, m = 1,2,...,M$.

\end{theorem}

\begin{proof}

According to the definition of the Frobenius norm and $tr(AB) = tr(BA)$ for every matrix $A, B$, we have

\begin{equation} \label{eq3.12}
\begin{split}
\sum_{m=1}^M \: \Vert \widetilde{\mathbb{A}}_m - \mathbb{U} \mathbb{V}_m^T \Vert_F^2 \: &= \: \sum_{m=1}^M \: tr((\widetilde{\mathbb{A}}_m - \mathbb{U} \mathbb{V}_m^T) (\widetilde{\mathbb{A}}_m - \mathbb{U} \mathbb{V}_m^T)^T) \\
&= \: \sum_{m=1}^M \: tr(\widetilde{\mathbb{A}}_m \widetilde{\mathbb{A}}_m^T) + \sum_{m=1}^M \: tr(\mathbb{V}_m \mathbb{V}_m^T) \\
& \quad \quad - 2 \sum_{m=1}^M \: tr(\widetilde{\mathbb{A}}_m^T \mathbb{U} \mathbb{V}_m^T)
\quad (since \: \mathbb{U} \: has \: orthonormal \: columns). 
\nonumber
\end{split}
\end{equation}

\noindent Obviously the first term on the right hand side of the equations above is a constant. Therefore, the minimization problem in Eq. (\ref{eq3.11}) is equivalent to solve

\begin{equation} 
\mathop{\min}\limits_{\mathbb{U}, \mathbb{V}_m} \: \sum_{m=1}^M \: tr(\mathbb{V}_m \mathbb{V}_m^T) \: - \: 2 \sum_{m=1}^M \: tr(\widetilde{\mathbb{A}}_m^T \mathbb{U} \mathbb{V}_m^T). 
\nonumber
\end{equation}

\noindent By the property of the trace of matrices, the minimizers of Eq. (\ref{eq3.11}) satisfy that $\mathbb{V}_m = \widetilde{\mathbb{A}}_m^T \mathbb{U}$ and $\mathbb{V}_m^T = \mathbb{U}^T \widetilde{\mathbb{A}}_m$, for any $m = 1,2,...,M$.

\end{proof}

Since $\mathbb{V}_m$ is uniquely determined by $\mathbb{U}$ and $\widetilde{\mathbb{A}}_m$ according to Theorem \ref{th3.3},  our next target is to obtain the semi-orthogonal matrix $\mathbb{U}$. 

\begin{theorem}\label{th3.4}

Let $\mathbb{U}$ be the optimal solution to the minimization problem in Eq. (\ref{eq3.11}), then it consists of the eigenvectors of the matrix

\begin{equation} \label{eq3.12}
\mathbb{N} = \sum_{m=1}^M \widetilde{\mathbb{A}}_m \widetilde{\mathbb{A}}_m^T,
\end{equation}

\noindent corresponding to the first $k$ largest eigenvalues.
\end{theorem}

\begin{proof}

The previous minimization problem in in Eq. (\ref{eq3.11}) is equivalent to the following optimization problem according to Theorem \ref{th3.3}:

\begin{equation}
\mathop{\max}\limits_{\mathbb{U}} \: \sum_{m=1}^M \:tr( \mathbb{U} \mathbb{U}^T \widetilde{\mathbb{A}}_m \widetilde{\mathbb{A}}_m^T). 
\nonumber
\end{equation}

\noindent By the properties of the trace of matrices, we have that

\begin{equation}
\begin{split}
\sum_{m=1}^M \:tr( \mathbb{U} \mathbb{U}^T \widetilde{\mathbb{A}}_m \widetilde{\mathbb{A}}_m^T) \: 
&= \: \sum_{m=1}^M \:tr( \mathbb{U}^T \widetilde{\mathbb{A}}_m \widetilde{\mathbb{A}}_m^T \mathbb{U}) \\
&= \: tr( \mathbb{U}^T (\sum_{m=1}^M \widetilde{\mathbb{A}}_m \widetilde{\mathbb{A}}_m^T) \mathbb{U}) \\
&= \: tr( \mathbb{U}^T \mathbb{N} \mathbb{U}).
\nonumber
\end{split}
\end{equation}

\noindent Since $\mathbb{N}$ is symmetric, we conclude that the maximum of the optimization problem above is achieved only if $\mathbb{U}$ is chosen to have the first $k$
eigenvectors of $\mathbb{N}$.

\end{proof}

We summarize the results above for computing the matrices $\mathbb{U}$ and $\{\mathbb{V}_m\}_{m=1}^M$ in Algorithm \ref{alg2}.

\begin{algorithm}[h]
\caption{SVD-based Low Rank Approximation Method of Matrices} \label{alg2}
\begin{algorithmic}[1]
	
\Require
Matrices $\{\widetilde{\mathbb{A}}_m\}_{m=1}^M$, and the dimension reduction ratio $\tau$
\Ensure
Matrices $\mathbb{U}$ and $\{\mathbb{V}\}_{m=1}^M$

\State Determine the rank $k = \lceil \tau N \rceil$. 

\State Form matrix $\mathbb{N}$ as Eq. (\ref{eq3.12}).

\State Compute $\mathbb{U}$ consisting of the $k$ eigenvectors of $\mathbb{V}$ corresponding to the first $k$ largest eigenvalues.

\For{$m = 1,...,M$}    

\State $\mathbb{V} \: = \: \widetilde{\mathbb{A}}_m^T \; \mathbb{U}$.

\EndFor

\State \Return matrices $\mathbb{U}$ and $\{\mathbb{V}\}_{m=1}^M$

\end{algorithmic}
\end{algorithm}

In Algorithm \ref{alg2}, we find that the formulation of the matrix $N$ in Line 2 and $\{\mathbb{V}\}_{m=1}^M$ in Lines 4-6, which take the computational time of $O(MN^2)$ and $O(Mk^2)$ respectively. Hence, the total computational cost is $O(M(N^2+k^2))$, as compared with SVD consuming time of $O(MN^4)$. Meanwhile, Line 3 in Algorithm \ref{alg1} involves an eigenvalue problem with the size of $N^2$, while the size is $MN^2$ in the traditional SVD method. A summary of the computational and space complexity of the two algorithms is given in Table \ref{tab1}.

\begin{table}[!ht] 
\setlength{\abovecaptionskip}{0pt}
\setlength{\belowcaptionskip}{10pt}
\centering 
\caption{Comparison of SVD, Algorithm \ref{alg1}, and Algorithm \ref{alg2}: $M, N, k, I$ denote the amount of samples, dimensions of the previous and low rank data, and the number of iterations in the while loop from Line 2 to Line 6 in Algorithm \ref{alg1} respectively.}
\begin{tabular}{ccc} \hline 
Method & Time & Space  \\ \hline
Traditional SVD & $O(MN^4)$ & $MN^2$ \\
Algorithm \ref{alg1} & $O(M(2IN^2+k^2))$ & $2Nk + Mk^2$ \\
Algorithm \ref{alg2} & $O(M(N^2+k^2))$ & $Nk + MNk$ \\ \hline
\end{tabular}\label{tab1}
\end{table}

\subsubsection{LRAM-based Method for Perturbed Linear Systems from ellipic SPDEs}

In Algorithm \ref{alg3}, we present the pseudo-code for solving the random linear system with perturbations in the coefficient matrices in Eq. (\ref{eq3.3}) formed by the probabilistic discretization of SPDE, based on low rank matrix approximation technique.  

\begin{algorithm}[h]
\caption{LRAM-based Method for Perturbed Linear Systems} 
\label{alg3}
\begin{algorithmic} [1]
	
\Require
Original stiffness matrix $\overline{\mathbb{A}}$, the perturbed stiffness matrix $\widetilde{\mathbb{A}}(\omega)$, the load vector $\bm{b}$, the amount of MC realizations $M$, and dimension reduction ratio $\tau$.
\Ensure
Approximation of QoI $\mu$.

\State Compute the original numerical solution $\overline{\bm{u}} = {\overline{\mathbb{A}}}^{-1} \bm{b}$ and sample the perturbed stiffness matrix $\widetilde{\mathbb{A}}_m, m = 1,...,M$.

\State Determine the rank $k = \lceil \tau N \rceil$ and make the optimal rank-$k$ approximation $\mathbb{U}, \{\mathbb{W}_m\}_{m=1}^M$ to matrices $\{\widetilde{\mathbb{A}}_m\}_{m=1}^M$ by Algorithm \ref{alg2}.

\For{$m = 1,...,M$}    

\State Compute the sample solutions $\bm{u}_{approx}^m$ by Eq. (\ref{eq3.6}).

\EndFor

\State \Return the QoI estimation $\mu \; = \; \sum_{m=1}^M \bm{u}_{approx}^m$.

\end{algorithmic}
\end{algorithm}

As is stated before, the high computational cost and memory of the traditional direct method to solve the perturbed linear system are significantly reducedin Algorithm \ref{alg3}. In Algorithm \ref{alg3}, we use the smaller sketches $\mathbb{U}$ and $\{\mathbb{W}_m\}_{m=1}^M$ with the determination of $\widetilde{\mathbb{A}}_m = \mathbb{U} \mathbb{W}_m$ to replace the $M$ perturbed stiffness matrices. It requires just $Nk + MNk$ scalars to store the rank-$k$ matrices. Actually, the single matrix $\mathbb{N}$ recovers enough information we desired from $\{\widetilde{\mathbb{A}}_m\}_{m=1}^M$, which also leads to the low memory storage requirement. The compression ratio $r$ reads

\begin{equation} \label{eq3.13}
\begin{split}
r \: & \:= \frac{Nk + MNk}{MNN} \:= \frac{k}{N} \: (1 + \frac{1}{M})\\
     & \:= \tau \: (1 + \frac{1}{M}) \enspace \rightarrow \: \tau, \quad as \enspace M \rightarrow \infty. 
\end{split}
\end{equation}

Both the speed-up of the matrix computation and the data compression ratio $r$ in Eq. (\ref{eq3.13}) level up as the dimension reduction ratio $\tau$ decreases. However, too small value of $\tau$ may lead to loss of information intrinsic in the original matrices. We will discuss this trade-off in the following sections.

Once the perturbation $\widetilde{\mathbb{A}}_m$ is approximated by $\widetilde{\mathbb{A}}_m= \mathbb{U} \mathbb{V}_m^T$, the perturbed linear system in Eq. (\ref{eq3.5}) is transformed into the following equations, where we denote $\mathbb{W}_m = \mathbb{V}_m^T$ for notational simplicity,

\begin{equation} \label{eq3.14}
(\overline{\mathbb{A}} \: + \: \mathbb{U} \mathbb{W}_m) \: \bm{u}_m \: = \: \bm{b}, \enspace m = 1,...,M.
\end{equation}

The final obstacle is to obtain the inverses of the  large-scale sparse matrices $(\overline{\mathbb{A}} + \mathbb{U} \mathbb{W}_m)$. Inspired by the Shermann-Morrison-Woodbury formula [\citealp{sherman1950adjustment, woodbury1950inverting, bartlett1951inverse}] , we approximate the solution to of Eq. (\ref{eq3.14}) 

$$ \bm{u}_m \: = \: (\overline{\mathbb{A}} \: + \: \mathbb{U} \mathbb{W}_m)^{-1} \bm{b},\\$$

by

\begin{equation} \label{eq3.15}
\begin{aligned}
& \bm{u}_{approx}^m \: = \: [\overline{\mathbb{A}}^{-1} - \overline{\mathbb{A}}^{-1} \mathbb{U} (\mathbb{I}_k + \mathbb{W}_m \overline{\mathbb{A}}^{-1} \mathbb{U})^{-1} \mathbb{W}_m \overline{\mathbb{A}}^{-1}] \bm{b},\\
& \qquad \quad \enspace = [\mathbb{I}_N - \overline{\mathbb{A}}^{-1} \mathbb{U} (\mathbb{I}_k + \mathbb{W}_m \overline{\mathbb{A}}^{-1} \mathbb{U})^{-1} \mathbb{W}_m] \overline{\mathbb{A}}^{-1} \bm{b}, \\
& \qquad \quad \enspace =  \overline{\bm{u}} - \overline{\mathbb{A}}^{-1} \mathbb{U} \mathbb{Y}_m \mathbb{W}_m \overline{\bm{u}}, \qquad \qquad m = 1,..., M,
\end{aligned}
\end{equation}

\noindent where  $\mathbb{Y}_m := (\mathbb{I}_k + \mathbb{W}_m \overline{\mathbb{A}}^{-1} \mathbb{U})^{-1}$ and  $\overline{\bm{u}} = \overline{\mathbb{A}}^{-1} \bm{b}$ being the unperturbed numerical solution. 

Note that the final numerical solution $\bm{u}_{approx}^m$ in  Eq. (\ref{eq3.15}) requires only $M$ matrix inversions of dimension $k$, which is significantly smaller than the matrix size $N$.  Moreover, Eq. (\ref{eq3.15}) indicates that the numerical solution $\bm{u}_{approx}^m$ can be viewed as a perturbed formulation of $\overline{\bm{u}}$, which is fixed. 

\begin{remark}
It is also feasible to obtain the inverses of the sparse large-scale matrices $(\overline{\mathbb{A}} + \mathbb{U} \mathbb{W}_m)^{-1}$ based on the properties of matrix series, i.e., if the spectral radius $\rho(\widetilde{\mathbb{A}}_m \overline{\mathbb{A}}^{-1}) < 1$, then the following series representation converges, which is majorized by

\begin{equation} 
(\mathbb{I} \: + \: \overline{\mathbb{A}}^{-1} \widetilde{\mathbb{A}}_m)^{-1} \enspace = \enspace \sum_{k=0}^\infty \; \left( - \; \overline{\mathbb{A}}^{-1} \widetilde{\mathbb{A}}_m \right)^k.
\nonumber
\end{equation}

The numerical solution in  Eq. (\ref{eq3.15}) is
\begin{equation} \label{eq3.15}
\begin{aligned}
\bm{u}_m \enspace & = \enspace (\overline{\mathbb{A}} \: + \: \widetilde{\mathbb{A}}_m)^{-1} \bm{b} \enspace = \enspace \left( \overline{\mathbb{A}}(\mathbb{I} \: + \: \overline{\mathbb{A}}^{-1} \widetilde{\mathbb{A}}_m) \right)^{-1} \bm{b}\\
& = \enspace (\mathbb{I} + \overline{\mathbb{A}}^{-1} \widetilde{\mathbb{A}}_m)^{-1} \overline{\mathbb{A}}^{-1} \bm{b} \enspace = \enspace (\mathbb{I} + \overline{\mathbb{A}}^{-1} \widetilde{\mathbb{A}}_m)^{-1} \; \overline{\bm{u}} \\
& = \enspace \sum_{k=0}^\infty ( -\; \overline{\mathbb{A}}^{-1} \widetilde{\mathbb{A}}_m)^k \; \overline{\bm{u}} \enspace = \enspace \sum_{k=0}^\infty ( - \; \overline{\mathbb{A}}^{-1} \mathbb{U} \mathbb{W}_m)^k \; \overline{\bm{u}} \\
& \approx \enspace \sum_{k=0}^K ( - \; \overline{\mathbb{A}}^{-1} \mathbb{U} \mathbb{W}_m)^k \; \overline{\bm{u}}, \qquad \qquad m = 1,..., M,
\end{aligned}
\end{equation}
where $K$ is order of truncation. 

The convergence of the series representation is guaranteed by the sufficiently small $\Vert \widetilde{\mathbb{A}}_m \Vert$ due to the virtue of $\rho(\widetilde{\mathbb{A}}_m \overline{\mathbb{A}}^{-1}) = \Vert \widetilde{\mathbb{A}}_m \overline{\mathbb{A}}^{-1} \Vert \le \Vert \widetilde{\mathbb{A}}_m \Vert \Vert\overline{\mathbb{A}}^{-1} \Vert < 1$. Such series-based approach is well suited for solving stochastic linear systems with perturbations in the coefficient matrices, since it generally implies that tiny perturbations are employed to these coefficients and thus guarantee the sufficiently small spectral radius, i.e. $\rho(\widetilde{\mathbb{A}}_m \overline{\mathbb{A}}^{-1}) << 1$. However, this approach requires a sufficiently large truncation coefficient $K$ to ensure the accuracy, while it also leads to high computational expense.

\end{remark}

\section{Applications on Elliptic SPDE and SOCP}

In this section, we consider two specific applications: the elliptic PDE with random diffusion coefficient and the stochastic optimal control problem (SOCP) governed by the elliptic SPDE.

\subsection{Stochastic Elliptic Partial Differential Equations} 

In a bounded, Lipschitz domain $D \in \mathbb{R}^d, d = 1,2,3$, we consider the following elliptic PDE with random diffusion coefficient: find a stochastic function $u : D \times \Omega \rightarrow \mathbb{R}$ such that the following equation holds with the homogeneous Dirichlet boundary condition:

\begin{equation}\label{eq4.1}
\left\{
\begin{aligned}
 - \nabla (a(\bm{x},\omega)\nabla u(\bm{x},\omega)) = f(\bm{x}), & \quad \bm{x} \in D, \omega \in \Omega,\\
 u(\bm{x},\omega) = 0, & \quad \bm{x} \in \partial D, \omega \in \Omega,
\end{aligned}
\right.
\end{equation}

\noindent where the diffusion coefficient $a$ is an almost surely continuous and positive random field on $D$, and we assume for simplicity that the force $f \in L^2(D)$.

We further assume the diffusion coefficient has the following form
\begin{equation} \label{eq4.2}
a(\bm{x}, \omega) = \overline{a}(\bm{x}) + \widetilde{a}(\bm{x}, \omega).
\end{equation}

\noindent where $a(\cdot, \omega)$ consists of its expectation $\overline{a}(\bm{x})>  0$ and a random field $\widetilde{a}(\bm{x}, \omega) \in \mathcal{L}^2(D)$ standing for its perturbation. To ensure the existence and uniqueness of the problem in Eq. (\ref{eq4.1}), we assume the stochastic diffusion coefficient $a(\cdot, \omega)$  satisfies the following assumptions.

\begin{assumption} [Regularity of Coefficients] \label{ass4.1}

There exists constants $a_{\min}$ and $a_{\max}$ such that the stochastic diffusion coefficient $a$ is uniformly elliptic, i.e.

\begin{equation}
0 \, < \, a_{\min} \leq a(\bm{x}, \omega) \leq a_{\max} \, < \, \infty, \quad a.e. \enspace (\bm{x}, \omega) \enspace in \enspace (D, \Omega).
\nonumber
\end{equation}
  
\end{assumption} 

\begin{assumption} \label{ass4.2}

The domain $D \subset \mathbb{R}^d$ is polygonal convex and the random field $a$ satisfies $a(\cdot, \omega) \in H_0^1(D)$ for a.e. $\omega \in \Omega$ with $ess \sup_{\omega} \Vert a(\cdot, \omega) \Vert_{H_0^1(D)} \, < \, \infty.$

\end{assumption}

Then we recall the well posedness of the stochastic elliptic PDE \citep{babuska2004galerkin, lord2014introduction}.

\begin{lemma} [Well Posedness of Eq. (\ref{eq4.1})] \label{le4.3}

Let Assumption \ref{ass4.1} and \ref{ass4.2} hold. If $f \in \mathcal{L}^2(D)$, then the equation admits a unique and bounded solution $u \in \mathcal{L}^2(D)$. There also exists a constant $C$, independent of $f$, s.t.

\begin{equation}
\Vert u \Vert_{\mathcal{L}^2(D)} \leq C \Vert f \Vert_{\mathcal{L}^2(D)}.
\nonumber
\end{equation}

\end{lemma} \label{le2}

For notational simplicity, we introduce the weak formulation of the stochastic elliptic PDE in Eq. (\ref{eq4.1}): find $u \in \mathcal{H}_0^1(D)$ satisfying that

\begin{equation} \label{eq4.3}
b[u,v] = [f,v], \quad \forall v \in \mathcal{H}_0^1(D), 
\quad for \enspace a.e. \enspace \omega \in \Omega,
\end{equation}

\noindent where the bilinear forms are given by

\begin{equation} 
b[u,v] := \mathbb{E} \left[ \int_D a(\bm{x}, \omega) \nabla u(\bm{x}, \omega) \cdot \nabla v(\bm{x}) d\bm{x} \right],
\nonumber
\end{equation}

\noindent and

\begin{equation} 
[f,v] := \int_D f(\bm{x}) v(\bm{x}) d\bm{x}.
\nonumber
\end{equation}

By Assumption \ref{ass4.1}, \ref{ass4.2} and Lax-Milgram Lemma \citep{lax2016ix}, the existence and uniqueness of weak solution to Eq. (\ref{eq4.1}) \citep{lord2014introduction, evans2022partial} can be obtained.

\begin{lemma} [Existence and Uniqueness of Eq. (\ref{eq4.1})] \label{le4.4}

Let $f \in L^2(D)$, then there exists a unique solution to the weak formulation in Eq. (\ref{eq4.1}) in $\mathcal{H}_0^1(D)$.

\end{lemma} 

\subsubsection{FE Discretization}

First we introduce the discrete formulation of the stochastic elliptic PDE in Eq. (\ref{eq4.1}). The Monte Carlo finite element method (MCFEM) \citep{gunzburger2014stochastic} is adopted in the article to alleviate the curse of dimensionality. Specifically, we approximate the integral $\mathbb{E}[u]$ in Eq. (\ref{eq3.2}) numerically by sample averages of realizations corresponding to the independent identically distributed (i.i.d.) random inputs. On the other hand, the standard finite element method is used in spatial discretization with respect to $\bm{x} \in D$.

Let $V_h \in H_0^1(D)$ be the finite element spaces associated with a regular shape mesh $\mathcal{T}_h$, and let $\{\phi_j\}_{j=1}^N$ denote the basis functions. Then we have

\begin{equation} \label{eq4.4}
u_h (\bm{x},\omega) = \sum_{j=1}^N u_{j,h}(\omega) \phi_j(\bm{x}) \quad and \quad f_h (\bm{x}) = \sum_{j=1}^N f_{j,h} \phi_j(\bm{x}),
\end{equation}

\noindent and the weak formulation in Eq. (\ref{eq4.3}) reduces to

\begin{equation} \label{eq4.5}
\begin{aligned}
\sum_{j=1}^N u_{j,h}(\omega) \int_D \overline{a}(\bm{x}) \nabla \phi_j(\bm{x}) & \nabla \phi_i(\bm{x}) d\bm{x} \: + \: \sum_{j=1}^N u_{j,h}(\omega) \int_D \widetilde{a}(\bm{x}, \omega) \nabla \phi_j(\bm{x}) \nabla \phi_i(\bm{x}) d\bm{x} \\
& = \quad \sum_{j=1}^N f_{j,h} \int_D \phi_j(\bm{x}) \phi_i(\bm{x}) d\bm{x},
\end{aligned}
\end{equation}

\noindent for $i = 1,..., N$, and $N$ denotes the amount of FE basis.

Given the i.i.d sample realizations $\widetilde{a}^m(\bm{x}) := \widetilde{a}(\bm{x}, \omega_m)$ of the random field $\widetilde{a}(\cdot, \omega)$ in the diffusion coefficient for $m = 1,...,M$, we could obtain the i.i.d samples $u_h^m(\bm{x}) := u_h(\bm{x}, \omega_m)$ of the FE solutions by solving the following $M$ variational problems: 

\begin{equation} \label{eq4.6}
\begin{aligned}
\sum_{j=1}^N u_j^m \int_D \overline{a}(\bm{x}) & \nabla \phi_j(\bm{x}) \nabla \phi_i(\bm{x}) d\bm{x} + \sum_{j=1}^N u_j^m \int_D \widetilde{a}^m (\bm{x}) \nabla \phi_j(\bm{x}) \nabla \phi_i(\bm{x}) d\bm{x} \\
& = \sum_{j=1}^N f_j \int_D \phi_j(\bm{x}) \phi_i(\bm{x}) d\bm{x}, \quad i = 1,..., N, \enspace m = 1,..., M.
\end{aligned}
\end{equation}

For computational simplicity, the discrete formulation may be written as a linear system of algebraic equations:

\begin{equation} \label{eq4.7}
\mathbb{A}^m \bm{u}_h^m = \bm{b}, \quad m = 1,..., M,
\end{equation}

\noindent where the stiffness matrix $\mathbb{A}^m, \overline{\mathbb{A}}, \widetilde{\mathbb{A}}^m \in \mathbb{R}^{N \times N}$, the mass matrix $\Phi \in \mathbb{R}^{N \times N}$ and the load vector $\bm{b} \in \mathbb{R}^N$ are respectively defined by

\begin{equation} \label{eq4.8}
\begin{aligned}
& \overline{\mathbb{A}}_{ij} = \int_D \overline{a}(\bm{x}) \nabla \phi_j(\bm{x}) \nabla \phi_i(\bm{x}) d\bm{x}, \quad \widetilde{\mathbb{A}}^m_{ij} \int_D \widetilde{a}^m (\bm{x}) \nabla \phi_j(\bm{x}) \nabla \phi_i(\bm{x}) d\bm{x},\\
& \mathbb{\mathbb{A}}^m = \overline{\mathbb{A}} + \widetilde{\mathbb{A}}^m, \quad \Phi_{ij} = \int_D \phi_j(\bm{x}) \phi_i(\bm{x}) d\bm{x}, \quad \bm{b} = \Phi\bm{f},\\
& \bm{u}_h^m = (u_{1,h}^m, ..., u_{N,h}^m)^T \quad and \quad \bm{f} = (f_{1,h},...,f_{N,h})^T.
\end{aligned}
\end{equation}

\subsubsection{LRAM-based Method to Stochastic Elliptic PDE}

Once the equations are fully discretized and the boundary conditions are imposed, the stiffness matrix $\overline{\mathbb{A}}$ has a full rank, while $\{\widetilde{\mathbb{A}}^m\}_{m=1}^M$ are rank-deficient matrices due to the existence of finite element boundary points. We make use of the low-rank structures of the collection of matrices $\{\widetilde{\mathbb{A}}^m\}_{m=1}^M$ to guarantee the high accuracy of the algorithm and solve the linear system in Eq. (\ref{eq4.7}) by Algorithm \ref{alg3}. Finally, the QoI is estimated by

\begin{equation} \label{eq4.9}
\mu_{h,M,\tau}(\bm{x}) := \frac{1}{M} \sum_{m=1}^M \bm{u}_{h,\tau}^m(\bm{x}) = \frac{1}{M} \sum_{m=1}^M \sum_{j=1}^N u_{j,h,\tau}^m \phi_j(\bm{x}),
\end{equation}

\noindent and the linear system in Eq. (\ref{eq4.7}) is solved as

\begin{equation} \label{eq4.10}
\bm{u}_{h,\tau}^m \: = \: [\mathbb{I}_N - \overline{\mathbb{A}}^{-1} \mathbb{U} \mathbb{Y}_m \mathbb{W}_m] \overline{\bm{u}}, \quad m = 1,..., M.
\end{equation}

Therefore, we derive an efficient and fast numerical solution in Eq. (\ref{eq4.10}) for the elliptic SPDE problem, since the calculations of $N \times N$ inverse matrices are reduced to $k \times k$ inversion. We also require less storage for large-scale matrices due to the low rank matrix approximation method in Algorithm \ref{alg2}.

We summarize the procedure to the pseudo-code for solving the stochastic elliptic partial differential equations in Eq. (\ref{eq4.1}) based on low rank matrix approximation technique in Algorithm \ref{alg4}.

\begin{algorithm}[h]
\caption{LRAM-based Method for the Elliptic SPDE in Eq. (\ref{eq4.1})} 
\label{alg4}
\begin{algorithmic} [1]
	
\Require
A tessellation $\mathcal{T}_h$ of $D$, stochastic diffusion coefficient $a(\bm{x},\omega)$, the force $f(\bm{x})$, the amount of MC realizations $M$, and dimension reduction ratio $\tau$.
\Ensure
Approximation of QoI $\mu_{h,M,\tau}$.

\State Construct $V_h \in H_0^1(D)$ and sample the random field in diffusion coefficient $\widetilde{a}^m, m = 1,...,M$.

\State Assemble the stiffness matrices $\overline{\mathbb{A}}, \{\widetilde{\mathbb{A}}^m\}_{m=1}^M$ and the load vector $\bm{b}$ by Eq. (\ref{eq4.8}).

\State Deal with the Dirichlet boundary conditions to $\overline{\mathbb{A}}$.

\State Determine the rank $k = \lceil \tau N \rceil$ and compute the sample solutions $\bm{u}_{h,\tau}^m$ by Algorithm \ref{alg3}.

\State \Return the estimation $\mu_{h,M,\tau}$ of the expectation $\mathbb{E}[u]$ by Eq. (\ref{eq4.9}).

\end{algorithmic}
\end{algorithm}

\subsubsection{Error Analysis}
In this section provides, we estimate the total approximation error of the QoI $\mathbb{E}[u]$ for the elliptic SPDE in Eq. (\ref{eq4.1}) by separately studying the errors due to the FE discretization, MC sampling and low rank matrix approximation, respectively.

\begin{theorem}[Spatial Discretization Error \citep{lord2014introduction}] \label{th4.5}

Let $V_h \in H_0^1(D)$ denote a piecewise linear finite element space defined by a regular shape mesh $\mathcal{T}_h$ and denote the solutions of Eq. (\ref{eq4.1}) and Eq. (\ref{eq4.5}) by $u, u_h$ respectively, then we have

\begin{equation} \label{eq4.11}
\Vert \mathbb{E}[u] - \mathbb{E}[u_h] \Vert_{\mathcal{L}^2(D)} \le C_1 h \Vert f \Vert_{\mathcal{L}^2(D)}.
\end{equation}

\end{theorem}

\begin{theorem}[Probability Discretization Error \citep{james1980monte}] \label{th4.6}

Let $\mu_{h,M} = \frac{1}{M} \sum_{m=1}^M u_h^m(\bm{x})$ denote the MCFEM numerical solution, then the error estimate is obtained by

\begin{equation} \label{eq4.12}
\Vert \mathbb{E}[u_h] - \mu_{h,M} \Vert_{\mathcal{L}^2(D)} \le C_2 \frac{1}{\sqrt{M}}.
\end{equation}

\end{theorem}

To obtain the error estimate from the low rank matrix approximation by Algorithm \ref{alg2}, we first look at the following theorem \citep{wedin1973perturbation}:

\begin{theorem} \label{th4.7}
Let $\mathbb{A}, \mathbb{B} \in \mathbb{R}^{n \times n}$ and denote the matrix $\mathbb{C} := \mathbb{B} - \mathbb{A}$, then we have

\begin{equation} \label{eq4.13}
\Vert \mathbb{B}^{-1} - \mathbb{A}^{-1} \Vert \enspace \le \enspace \mu \; \max \{ \Vert \mathbb{A}^{-1} \Vert_2^2, \; \Vert \mathbb{B}^{-1} \Vert_2^2 \} \; \Vert \mathbb{C} \Vert,
\end{equation}

\noindent where $\mu$ is a constant and $\Vert \cdot \Vert$ denotes any norm.

\end{theorem}

Based on Theorem \ref{th4.7} concerning the matrix inversion, we now provide an upper bound for the error from approximating the low rank matrices. Here we provide the formula of the root mean square reconstruction error of Algorithm \ref{alg3}

\begin{equation} \label{eq4.14}
RMSRE(M) \: := \: \sqrt{\frac{1}{M} \sum_{m=1}^M \Vert \widetilde{\mathbb{\mathbb{A}}}_m - \mathbb{\mathbb{U}} \mathbb{W}_m \Vert_F^2}.
\end{equation}

\begin{theorem}[Low Rank Matrix Approximation Error] \label{th4.8}

Define $\mu_{h,M,\tau}$ as Eq. (\ref{eq4.10}) and we provide the error estimation as 

\begin{equation} \label{eq4.15}
\Vert \mu_{h,M} - \mu_{h,M,\tau} \Vert_{\mathcal{L}^2(D)} \le C_3 \; \frac{1}{\sqrt{M}} \; RMSRE(M).
\end{equation}

\end{theorem}

\begin{proof}

Denote $\mathbb{B}^m := \overline{\mathbb{A}} + \mathbb{U} \mathbb{W}^m$ for notational convenience. By Jensen's inequality and Theorem \ref{th4.7}, the left hand side (LHS) goes as

\begin{equation}
\begin{aligned}
\Vert \mu_{h,M} - \mu_{h,M,\tau} \Vert_{\mathcal{L}^2(D)}^2 \quad & = \quad \Vert \frac{1}{M} \: \sum_{m=1}^M (u_h^m(\bm{x})) - u_{h,\tau}^m(\bm{x}))\Vert_{\mathcal{L}^2(D)}^2 \\
& \le \quad \frac{1}{M^2} \: \sum_{m=1}^M \: \Vert u_h^m(\bm{x})) - u_{h,\tau}^m(\bm{x}) \Vert_{\mathcal{L}^2(D)}^2 \\
& = \quad \frac{1}{M^2} \: \sum_{m=1}^M \: \Vert \left( (\mathbb{A}^m)^{-1} - (\mathbb{B}^m)^{-1} \right) \; \bm{b} \; \phi(\bm{x}) \Vert_{\mathcal{L}^2(D)}^2\\
& \le \quad \frac{1}{M^2} \: \Vert \bm{b} \Vert^2_2 \: \Vert \phi(\bm{x}) \Vert_{\mathcal{L}^2(D)}^2 \: \sum_{m=1}^M \: \Vert (\mathbb{A}^m)^{-1} - (\mathbb{B}^m)^{-1} \Vert_2^2 \\
& \le \quad \frac{1}{M^2} \: \mu^2 \: \Vert \bm{b} \Vert_2^2 \: \Vert \phi(\bm{x}) \Vert_{\mathcal{L}^2(D)}^2 \: \mathop{\max}\limits_{m = 1,...,M} \{ \Vert (\mathbb{A}^m)^{-1} \Vert^4_2, \Vert (\mathbb{B}^m)^{-1} \Vert^4_2 \} \: \sum_{m=1}^M \: \Vert \mathbb{A}^m - \mathbb{B}^m \Vert_2^2 \\
& \le \quad \frac{1}{M^2} \: \mu^2 \: \Vert \bm{b} \Vert_2^2 \: \Vert \phi(\bm{x}) \Vert_{\mathcal{L}^2(D)}^2 \: \mathop{\max}\limits_{m = 1,...,M} \{ \Vert (\mathbb{A}^m)^{-1} \Vert^4_2, \Vert (\mathbb{B}^m)^{-1} \Vert^4_2 \} \: \sum_{m=1}^M \: \Vert \widetilde{\mathbb{A}}^m - \mathbb{U} \mathbb{W}^m \Vert_F^2, \\
\nonumber
\end{aligned}
\end{equation}

\noindent where $\phi(\bm{x})$ is the finite element basis function, and the derivation results from the properties of the norm of matrices $\Vert \mathbb{A}\mathbb{B} \Vert \le \Vert \mathbb{A} \Vert \Vert \mathbb{B} \Vert$, $\Vert \mathbb{A} + \mathbb{B} \Vert \le \Vert \mathbb{A} \Vert + \Vert \mathbb{B} \Vert$, and $\Vert \mathbb{A} \Vert_2^2 \le \Vert \mathbb{A} \Vert_F^2$.

\noindent Thus, we obtain  

\begin{equation}
\begin{aligned}
\Vert \mu_{h,M} - \mu_{h,M,\tau} \Vert_{\mathcal{L}^2(D)} \enspace \le \enspace \left( \Vert \mu_{h,M} - \mu_{h,M,\tau} \Vert_{\mathcal{L}^2(D)}^2 \right)^{\frac{1}{2}}
\enspace \le \enspace C_3 \: \frac{1}{\sqrt{M}} \: RMSRE(M), \\
\nonumber
\end{aligned}
\end{equation}

\noindent where $C_3 := \mu \; \Vert \bm{b} \Vert_2 \; \Vert \phi(\bm{x}) \Vert_{\mathcal{L}^2(D)} \; \mathop{\max}\limits_{m = 1,...,M} \{ \Vert (\mathbb{A}^m)^{-1} \Vert^2_2, \Vert (\mathbb{B}^m)^{-1} \Vert^2_2 \}$.

\end{proof}

Based on the estimations from Theorem \ref{th4.5}, \ref{th4.6} and \ref{th4.8}, we have the following error analysis estimation.

\begin{theorem}[Error Analysis] \label{th4.9}

For any $M \in \mathbb{N}$, $\tau \in (0,1)$ and $u \in \mathcal{L}^2(D)$ holds that

\begin{equation} \label{eq4.16} 
\Vert \mathbb{E}[u] - \mu_{h,M,\tau} \Vert_{\mathcal{L}^2(D)} \quad = \quad O(h) \enspace + \enspace O(\frac{1}{\sqrt{M}}) \enspace + \enspace O(RMSRE(M)).
\end{equation}

\end{theorem}

\begin{proof}

By the triangle inequality, then we have that that

\begin{equation}
\begin{aligned} 
\Vert \mathbb{E}[u] - \mu_{h,M,\tau} \Vert_{\mathcal{L}^2(D)} \quad & \le \quad
\Vert \mathbb{E}[u] - \mathbb{E}[u_h] \Vert_{\mathcal{L}^2(D)} \enspace + \enspace \Vert \mathbb{E}[u_h] - \mu_{h,M} \Vert_{\mathcal{L}^2(D)} \enspace + \enspace \Vert \mu_{h,M} - \mu_{h,M,\tau} \Vert_{\mathcal{L}^2(D)} \\
& \le \quad C_1 h \Vert f \Vert_{\mathcal{L}^2(D)} \enspace + \enspace C_2 \frac{1}{\sqrt{M}} \enspace + \enspace C_3 \: \frac{1}{\sqrt{M}} \: RMSRE(M) \\
& = \quad O(h) \enspace + \enspace O(\frac{1}{\sqrt{M}}) \enspace + \enspace O(RMSRE(M)),
\nonumber
\end{aligned}
\end{equation}

\noindent where the reconstruction error $RMSRE$ is related to the dimension reduction ratio $\tau$.

\end{proof}

\subsubsection{Numerical Experiments on the Elliptic SPDE}

Consider the two-dimensional stochastic elliptic boundary value problem as illustrated in Eq. (\ref{eq4.1}), we test the numerical example with the following settings: let the spatial variable $\bm{x} = [x,y]^T$, and give a constant source term $f(\bm{x}) = 1$ for $\bm{x} \in D$. Uncertainties of the system comes from the permeability field $a(\bm{x}, \omega)$, which has the form

\begin{equation} \label{eq4.17}
a(\bm{x}, \omega) \: = \: 1 + \epsilon \sigma(\bm{x}, \omega),
\end{equation}

\noindent where the magnitude of the perturbation $\epsilon  =  0.2$, and $\sigma(\bm{x}, \omega) \in \mathcal{L}^2(D)$ is a random process with each component following the standard normal distribution $N(0,1)$. Obviously, the stochastic coefficient $a(\cdot, \omega)$ satisfies Assumption \ref{ass4.1}.

\begin{figure}[ht]
\centering
\includegraphics[width=0.5\textwidth]{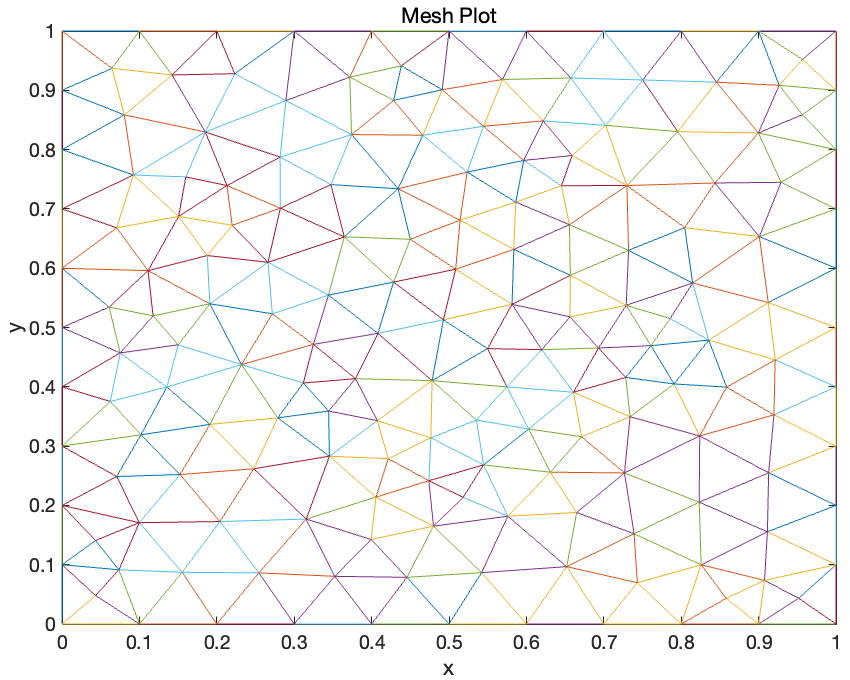}
\caption{Triangularized mesh plot of the domain $D = [0,1]^2$ .}
\label{fig1}
\end{figure}

Our goal is to compute the QoI defined in Eq. (\ref{eq4.9}). For the numerical implementation of Algorithm \ref{alg4}, we take the finite element mesh size $h = 0.1$, and the number of the finite element nodes $N = 665$, which is presented in the Figure \ref{fig1}, and set the amount of the MC realizations $M = 500$ in the MCFEM method. We carry
out simulations by using MATLAB R2022a software on an Apple M1 machine with 8GB of memory.

\begin{figure}[ht]
\centering
\includegraphics[width=1 \textwidth]{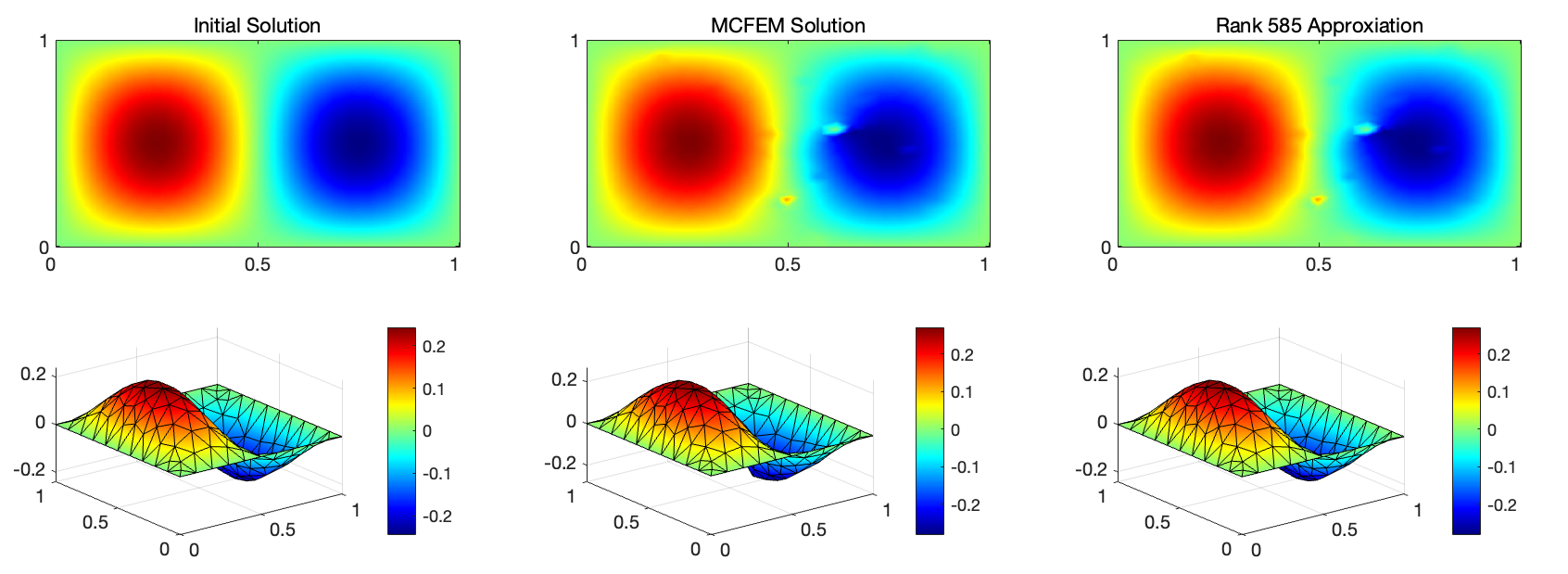}
\caption{The initial solution $\overline{\bm{u}}$ (left), the QoI obtained by MCFEM $\mu_{h,M}$ (middle), and the QoI obtained by Algorithm \ref{alg4} $\mu_{h,M,\tau}$ (right).}
\label{fig4.2}
\end{figure}

The numerical results of solving the elliptic SPDE in Eq. (\ref{eq4.1}) via Algorithm \ref{alg4} with the dimension reduction ratio $\tau = 0.88$ are demonstrated in Figure \ref{fig4.2}. Compared with the left one, the middle plot has a small gap in the margin. And when we select the dimension reduction ratio $\tau = 0.88$, i.e., each $665 \times 665$ sample of the stochastic stiffness matrices $\widetilde{\mathbb{A}}_m$ is factorized by two $665 \times 585$ low rank matrices $\mathbb{U}$ and $\mathbb{V}$ by $\widetilde{\mathbb{A}}_m = \mathbb{U} \mathbb{V}^T$, then we obtain the rank-$585$ approximation of $u$ visualized in the right plot. The middle and right plot are approximately identical and the error $\Vert \mu_{h,M} - \mu_{h,M,\tau} \Vert_2 = 7.92 \times 10^{-13}$, which demonstrates the effectiveness and validity of Algorithm \ref{alg4} in solving the stochastic partial differential equation.

\begin{figure}[ht]
\centering
\includegraphics[width=1 \textwidth]{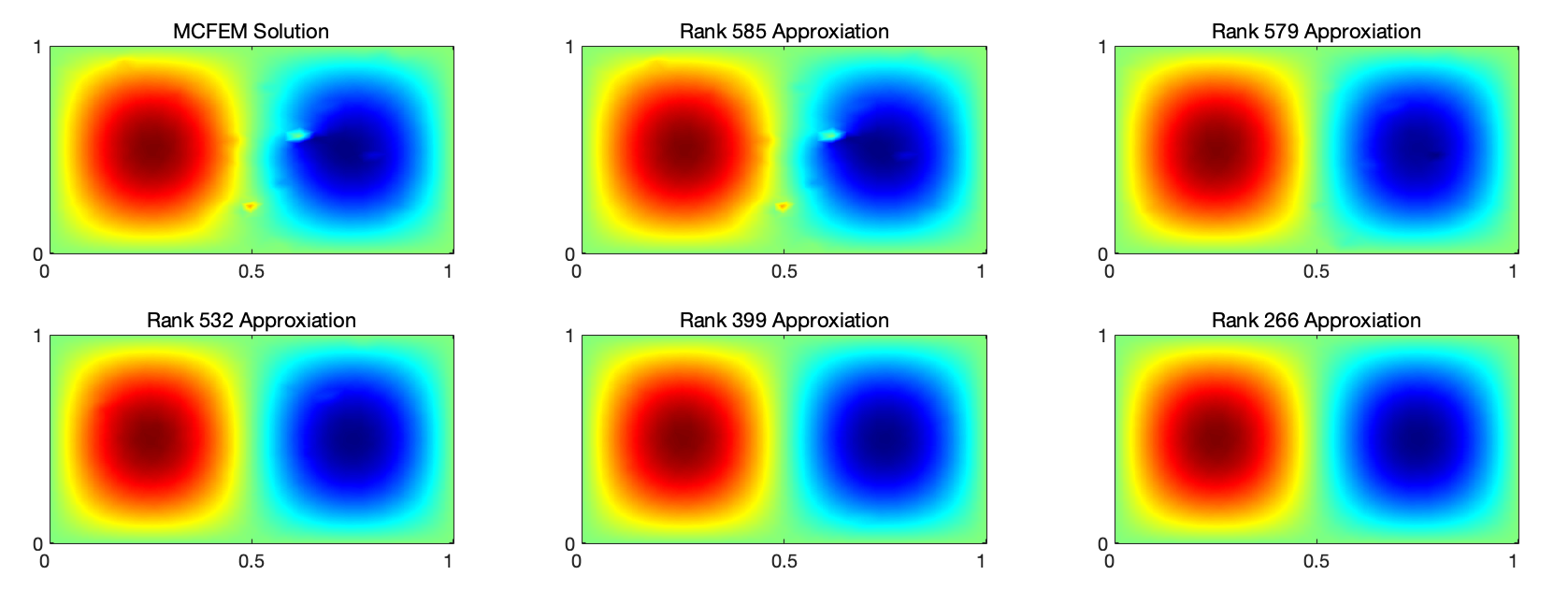}
\caption{The numerical QoI $\mu_{h,M,\tau}$ with 5 different dimension reduction ratio $\tau = 0.88$, $\tau = 0.87$, $\tau = 0.8$, $\tau = 0.6$, and $\tau = 0.4$.}
\label{fig4.3}
\end{figure}

We also examine the sensitivity of Algorithm \ref{alg4} to the choice of the dimension reduction ratio $\tau$. Figure \ref{fig4.3} depicts the numerical QoI $\mu_{h,M,\tau}$ and some simulation results with 5 different
$\tau = 0.88, 0.87, 0.8, 0.6, 0.4$. Compared to $\tau = 0.88$, the visualizations of $\mu_{h,M,\tau}$ with the other dimension reduction ratios do not capture the perturbation caused by the uncertainty in the diffusion coefficient and look relatively different from the MCFEM numerical solution. 

Table \ref{tab2} presents the errors between $\mu_{h,M,\tau}$ and the MCFEM numerical QoI $\mu_{h,M}$ and the CPU elapsed time under different $\tau$. It shows that Algorithm \ref{alg3} with $\tau = 0.88$ has much higher computational accuracy than other dimension reduction ratios, while using lower $\tau$ will significantly save the computational time and storage. Under these observations, we conclude two possible applications of Algorithm \ref{alg4}: 

\begin{itemize}

    \item Construct a high precision algorithm for solving stochastic SPDE. When we select a relatively high dimension reduction ratio, we will obtain the numerical results with a high computational accuracy. Simultaneously, the algorithm has relatively low requirement for CPU memory and computational complexity. We reduce the dimensions of the matrices $\{\widetilde{A}^m\}_{m=1}^M$ and simplify derivation for the numerical solution in Eq. (\ref{eq4.10}), in which we directly obtain the perturbed solution $\bm{u}_{h,M,\tau}$ by using the deterministic solution $\overline{\bm{u}}$ instead of solving a new linear system.

    \item Serve as a pre-processing for solving stochastic SPDE. In practical industrial implementation, it is quite difficult and time-consuming to compute a high-precision numerical solution. Instead, we first obtain an approximation which captures the main sketch of the exact solution and then carve it gradually. Algorithm \ref{alg4} with a relatively low dimension reduction ratio enjoys rapid computational speed and low storage reserves, since the dimensions of the matrices are significantly reduced by low rank matrix approximation technique in Algorithm \ref{alg2}.
\end{itemize}

\begin{table}[ht]
\centering 
\caption{Simulation results for Algorithm \ref{alg4} about the CPU elapsed time and the error. }
\begin{tabular}{cccccccc}
\hline
$\tau$ & $1.0$ & $0.88$ & $0.87$ & $0.8$ & $0.6$ & $0.4$ \\
Time (s) & $31.8702$ & $27.0772$ & $24.7071$ & $20.8188$ & $12.1475$ & $6.3270$ \\
Error & $1.2089 \times 10^{-12}$ & $1.2016 \times 10^{-12}$ & $0.5120$ & $0.5150$ & $0.5179$ & $0.5421$\\
\hline
\end{tabular}
\label{tab2}
\end{table}

It is also observed from Table \ref{tab2} that the error rises sharply as dimension reduction ratio $\tau$ falls below the critical point. The dramatic downward trend of the numerical efficiency results from the ill-conditioning of the stiffness matrix. The condition number serves as a measure of stability for linear systems \citep{rainer1998numerical}, and the perturbed stiffness matrices $\{\widetilde{\mathbb{A}}_m\}_{m=1}^M$ all have large condition numbers of $cond(\widetilde{\mathbb{A}}_m) \ge 3.5775 \times 10^{17}, m = 1, ..., M$ while we have $cond(\overline{\mathbb{A}}) = 339.3037$. The large value indicates the high sensitivity to uncertainties and perturbations in the linear system. Therefore, we will obtain a completely different solutions when there is a small change in stiffness matrices. 

The critical point of such two applications is $\tau = 0.88$ in this numerical settings. It is important to figure out how to determine the value of the critical point. Our idea is that it comes from the energy ratio of matrix $N = \sum_{m=1}^M \widetilde{\mathbb{A}}_m \widetilde{\mathbb{A}}_m^T$ 

\begin{equation}
e(k) \enspace := \enspace \sum_{i=1}^k \sigma_i^2 \: / \: \sum_{i=1}^N \sigma_i^2,
\nonumber
\end{equation}

\noindent where $\sigma_i$ denotes the $i$-th eigenvalue of the matrix $N$. The basis of Algorithm \ref{alg3} is the singular value decomposition, where the eigenvectors describe the directions of matrix transformation and the corresponding eigenvalues denote their importance. The energy ratio $e(k)$ guides us to determine an appropriate truncation index $k$ in the low rank matrix approximation, where we retain enough information from the original matrix $\widetilde{\mathbb{A}}(\omega)$ \citep{mcgivney2014svd}. In other words, we extract the main features of the matrices $\{\widetilde{\mathbb{A}}^m\}_{m=1}^M$ by taking the first $k = \lceil \tau N \rceil$ eigenvectors of $\mathbb{N}$ in Algorithm \ref{alg3}, since the first $k$ eigenvectors occupy considerably large energy of the matrix $\mathbb{N}$, i.e., they contain as much of information as we desire.

\begin{figure}[ht]
\centering
\includegraphics[width=0.95 \textwidth]{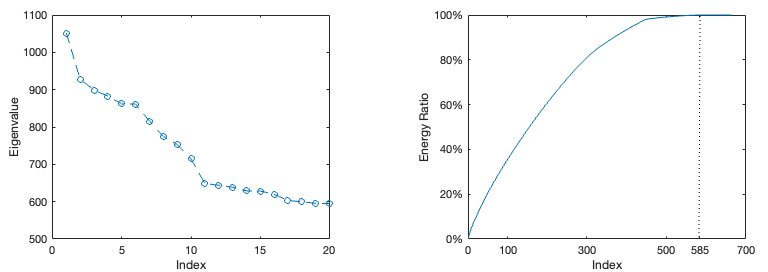}
\caption{The first 20 eigenvalues (left) and the energy ratios (right) of matrix $\mathbb{N}$.}
\label{fig4.4}
\end{figure}

The left plot in Figure \ref{fig4.4} depicts the first 20 eigenvalues of the matrix $\mathbb{N}$ listed in descending order, where the slopes of several domains are steep. The right plot demonstrates the energy ratios of the matrix $\mathbb{N}$. We observe that energy ratio $e(k)$ grows with the index $k$ and it reaches the top with $k = 585$. The simulation results reveal that the corresponding optimal rank-$k = 585$ approximation contains the whole information of the original stiffness matrix and thus we obtain the high-precision numerical solution by dimension reduction ratio $\tau = 585/665 \approx 0.88$. As a result, the value of the critical point $\tau$ could be determined by computing the eigenvalues and energy ratios of matrix $\mathbb{N}$.

\subsection{Stochastic Elliptic Control Problem}

In this application, we consider a stochastic optimal control problem, whose objective functional is of a velocity tracking type, and the governing system is defined by the elliptic PDE with a stochastic diffusion coefficient, i.e.,

\begin{equation} \label{eq4.18}
\min \mathcal{J}(u, f) \: := \: \mathbb{E} \left[\frac{1}{2} \Vert u - U \Vert^2  + \frac{\beta}{2} \Vert f \Vert^2 \right],
\end{equation}

\noindent subject to

\begin{equation} \label{eq4.19}
\left\{
\begin{aligned}
 - \nabla (a(\bm{x},\omega)\nabla u(\bm{x},\omega)) = f(\bm{x}), & \quad \bm{x} \in D, \omega \in \Omega,\\
 u(\bm{x},\omega) = 0, & \quad \bm{x} \in \partial D, \omega \in \Omega,
\end{aligned}
\right.
\end{equation}

\noindent where $u$ and $U$ denote the state variable and the deterministic desired state respectively, $f$ is a deterministic control, and $\beta$ is a small control penalty parameter. The nonempty admissible set $\mathcal{U}_{ad}$ is given by

\begin{equation}
\begin{aligned}
\mathcal{U}_{ad} \enspace &:= \enspace \mathcal{H}_0^1(D) \: \otimes \: \{ f \in L^2(D): f(\bm{x}) \ge 0, \forall x \in D \}.
\end{aligned}
\nonumber
\end{equation}

The weak form of the governing equation Eq. (\ref{eq4.19}) is given as Eq. \eqref{eq4.3}, and the existence and uniqueness of its weak solution are guaranteed by Lemma \ref{le4.3}. Therefor the weak formulation of the optimization problem  \eqref{eq4.18} is given as: determine the state $u^*$ and the control $f^*$ such that

\begin{equation} \label{eq4.20}
\mathop{\min}\limits_{(u,f) \: \in \: \mathcal{U}_{ad}}  \: \mathcal{J}(u,f), \qquad s.t. \quad b[u,v] = [f,v], \enspace \forall v \in H_0^1(D).
\end{equation}

Then the existence and uniqueness of the optimal solution of Eq. (\ref{eq4.18}) is present in the following lemma.

\begin{lemma} [\citealp{ali2017multilevel}] \label{le4.9}
Suppose that $\mathcal{U}_{ad}$ is nonempty. Then the distributed elliptic optimal control problem in Eq. (\ref{eq4.18}) - (\ref{eq4.19}) has a unique optimal solution $u^* \in \mathcal{U}_{ad}$ for almost surely $\omega \in \Omega$.
\end{lemma}

\subsubsection{Discretization of the Optimality System}

There are mainly two different strategies to solve the optimal control problems in literature: Discretize-then-Optimize approach \citep{leykekhman2012investigation, liu2019non, mathew2007analysis} and Optimize-then-Discretize approach \citep{luo2013some, neitzel2009strategies, tiesler2012stochastic}. Here we adopt the former strategy begin by discretizing our optimal control problem in Eq. (\ref{eq4.18}) and (\ref{eq4.19}) using the MCFEM method. In other words, for the stochastic process $u(\bm{x}, \omega)$, we will omit its explicit dependence on the probability space and treat each realization as a deterministic PDE. 

Given a collection of i.i.d. samples $\{\omega_m\}_{m=1}^M$ randomly drawn from $\Omega$ by the probability measure $\mathbb{P}$. Let $V_h \in H_0^1(D)$ denote the finite element subspace corresponding to the mesh size $h \; > \; 0$, then the fully discrete objective functional is stated as

\begin{equation} \label{eq4.21}
\min \widehat{\mathcal{J}}(u_{h,M}, f_h) \: := \: \frac{1}{M} \sum_{m=1}^M \left[\frac{1}{2} \Vert u_{h,M}(\omega_m) - U \Vert^2  + \frac{\beta}{2} \Vert f_h \Vert^2 \right],
\end{equation}

\noindent and the discretization of the weak formulation is given by

\begin{equation} \label{eq4.22}
b[u_{h,M}(\omega_m), v_h] = [f_h, v_h], \enspace \forall v_h \in V_h, \enspace m = 1,...,M.
\end{equation}

In short, the discrete optimal control problem is transformed to minimize the objective functional in Eq. (\ref{eq4.21}) subject to $M$ realizations of the discrete state system in Eq. (\ref{eq4.22}). As present in Section 4.1.2, we could further discretize the system by MCFEM method and perform dimensionality reduction towards the stochastic matrices $\{ \widetilde{\mathbb{A}}^m \}_{m=1}^M$. Then, by Algorithm \ref{alg3}, the state equation has the following algebraic form 

\begin{equation} \label{eq4.23}
\begin{aligned}
& (\overline{\mathbb{A}} + \mathbb{U} \mathbb{W}_m) \bm{u}_{h,\tau}^m = \Phi \bm{f}_h, \quad m = 1,..., M \\
\Rightarrow \quad & \bm{u}_{h,\tau}^m = (\overline{\mathbb{A}} + \mathbb{U} \mathbb{W}_m)^{-1} \Phi \bm{f}_h, \\
\Rightarrow \quad & \bm{u}_{h,\tau}^m = \left( \mathbb{I}_N - \overline{\mathbb{A}}^{-1} \mathbb{U} \mathbb{Y}_m \mathbb{W}_m \right) \overline{\mathbb{A}}^{-1} \Phi \bm{f}_h, \\
\Rightarrow \quad & \bm{u}_{h,\tau}^m = \mathbb{Z}_m \bm{f}_h,
\end{aligned}
\end{equation}

\noindent where we denote the  matrices

\begin{equation} \label{eq4.24}
\begin{aligned}
& \mathbb{Z}_m \: := \: \left(\mathbb{I}_N - \overline{\mathbb{A}}^{-1} \mathbb{U} \mathbb{Y}_m \mathbb{W}_m \right) \overline{\mathbb{A}}^{-1} \Phi, \\
& \mathbb{Y}_m \: := \: (\mathbb{I}_k + \mathbb{W}_m \overline{\mathbb{A}}^{-1} \mathbb{U})^{-1}. \\
\end{aligned}
\end{equation}
In Eq. (\ref{eq4.23}) and (\ref{eq4.24}), $\Phi, \overline{\mathbb{A}} \in \mathbb{R}^{N \times N}$ are the finite element mass matrix and the stiffness matrix corresponding to the deterministic diffusion coefficient $\overline{a}(\bm{x})$ respectively, $\mathbb{U}, \{ \mathbb{W}^m \}_{m=1}^M \in \mathbb{R}^{N \times \lceil \tau N \rceil}$ come from the low rank approximation of $\{ \widetilde{\mathbb{A}}^m \}_{m=1}^M$ with respect to $\widetilde{a}(\cdot, \omega)$, $\bm{f}_h = \left( f_{1,h}, ..., f_{N,h} \right)^T$ is the vector of $N$ degrees of freedom corresponding to the control $f_h (\bm{x})$ and the vector $\bm{u}_{h,\tau}^m = \left( u_{1,h,\tau}^m, ..., u_{N,h,\tau}^m \right)^T \in \mathbb{R}^N$ corresponds to the state $u_{h,M,\tau}^m (\bm{x})$ for $m = 1,...,M$, where the state and control are respectively given by

\begin{equation} \label{eq4.25}
u_{h,M,\tau}^m (\bm{x}) = \sum_{j=1}^N u_{j,h,\tau}^m \phi_j(\bm{x}) \quad and \quad f_h (\bm{x}) \: = \: \sum_{j=1}^N f_{j,h} \phi_j(\bm{x}). 
\end{equation}

It is shown in Eq. (\ref{eq4.23}) that the state $u$ depends on the control $f$, i.e. there exists a unique vector $\bm{u}_{h,\tau}^m$ for any $\bm{f}_h$. We denote $u_{h,M,\tau}^m = \mathcal{U}^m_\tau(f_h)$ to emphasize the dependence, where $\mathcal{U}^m_\tau$ is a linear operator. Therefore the discrete optimal control problem can be rewritten as

\begin{equation} 
\min \: \widehat{\mathcal{J}}(u_{h,M,\tau}^m, f_h), \quad s.t. \enspace u_{h,M,\tau}^m = \mathcal{U}^m_\tau(f_h).
\nonumber
\end{equation}

\subsubsection{LRAM-based Method for Stochastic Elliptic OCP}
We develop and analyze the fast algorithm in solving stochastic optimal control problems constrained by the elliptic SPDE. To reduce computational complexity, we plug the discrete constraints $u_{h,M,\tau}^m = \mathcal{U}^m_\tau(f_h)$ into objective functional and  obtain the reduced problem formulation as below

\begin{equation} \label{eq4.26}
\mathop{\min}\limits_{f_h} \: \widehat{\mathcal{J}}(\mathcal{U}^m_\tau(f_h), f_h).
\end{equation}

Such transformation allows an explicit elimination of constraints so that the stochastic optimal control problem becomes an unconstrained optimization problem, which has relatively low computational requirement and benefits from various kinds of gradient-based algorithms \citep{chen2022reduced} -- the first-order derivative method, such as the steepest descent method, and the second-order derivative method, like Newton's method, etc. 

By directly plugging the algebraic formulation in Eq. (\ref{eq4.23}), we obtain the following explicit formulation of the discretized and reduced objective functional

\begin{equation} \label{eq4.27}
\widehat{\mathcal{J}}(\bm{f}_h) \: = \: \frac{1}{M} \sum_{m=1}^M \frac{1}{2} (\bm{u}_{h,\tau}^m - \bm{U})^T \Phi (\bm{u}_{h,\tau}^m - \bm{U}) + \frac{\beta}{2} \bm{f}_h^T \Phi \bm{f}_h,
\end{equation}

\noindent where we define $\left( \bm{U} \right)_j = \int_D U(\bm{x}) \phi_j(\bm{x}) d\bm{x}$.

We first derive the gradient of the objective functional $\widehat{\mathcal{J}}(\bm{f}_h)$

\begin{equation} \label{eq4.28}
\nabla \widehat{\mathcal{J}}(\bm{f}_h) \: = \: \frac{\mathrm{d} \widehat{\mathcal{J}}}{\mathrm{d} \bm{f}_h} \: = \: \frac{1}{M} \sum_{m=1}^M \left( \frac{\partial \widehat{\mathcal{J}}}{\partial \bm{u}_{h,\tau}^m} \right)^T \frac{\mathrm{d} \bm{u}_{h,\tau}^m}{\mathrm{d} \bm{f}_h}  + \frac{\partial \widehat{\mathcal{J}}}{\partial \bm{f}_h},
\end{equation}

\noindent where by the direct differentiation of the formulations in Eq. (\ref{eq4.27}) and (\ref{eq4.23}), the partial derivatives $\frac{\partial \widehat{\mathcal{J}}}{\partial \bm{u}_{h,\tau}^m}$, $\frac{\partial \widehat{\mathcal{J}}}{\partial \bm{f}_h}$, and the sensitivity $\frac{\mathrm{d} \bm{u}_{h,\tau}^m}{\mathrm{d} \bm{f}_h}$ are easily determined as

\begin{equation} \label{eq4.29}
\frac{\partial \widehat{\mathcal{J}}}{\partial \bm{u}_{h,\tau}^m} \: = \: \Phi (Z^m \bm{f}_h - \bm{U}), \qquad \frac{\partial \widehat{\mathcal{J}}}{\partial \bm{f}_h} \: = \: \beta \Phi \bm{f}_h, \qquad \frac{\mathrm{d} \bm{u}_{h,\tau}^m}{\mathrm{d} \bm{f}_h} \: = \: Z^m.
\end{equation}
Substitute Eq. (\ref{eq4.29}) into (\ref{eq4.28}), the gradient of the objective functional $\widehat{\mathcal{J}}(\bm{f}_h)$ is deduced to

\begin{equation} \label{eq4.30}
\nabla \widehat{\mathcal{J}}(\bm{f}_h) \: = \: \frac{1}{M} \sum_{m=1}^M (Z^m)^T  \Phi (Z^m \bm{f}_h - \bm{U}) + \beta \Phi \bm{f}_h.
\end{equation}

Similarly, we can the derive the Hessian matrix of the objective functional $\widehat{\mathcal{J}}(\bm{f}_h)$ 

\begin{equation} \label{eq4.31}
\Delta \widehat{\mathcal{J}} \: = \: \frac{1}{M} \sum_{m=1}^M (Z^m)^T  \Phi Z^m + \beta \Phi.
\end{equation}

Therefore, we can update descend direction in the optimization loop based on the first- and second-order derivatives from Eq. (\ref{eq4.30}) and (\ref{eq4.31}). Note that the formulation of $\Delta \widehat{\mathcal{J}}$ implies that we don't need calculate the hessian matrix in each iteration since it is independent of $\bm{f}_h$. We present the pseudo-code of the gradient-based optimization algorithm to solve the stochastic optimal control problem governed by the elliptic SPDE in Algorithm \ref{alg5}.

\begin{breakablealgorithm}
\caption{LRAM-based Method for the stochastic optimal control problem in Eq. (\ref{eq4.18}) governed by the elliptic SPDE in Eq. (\ref{eq4.19}) via the gradient-based optimization algorithm and low rank matrix approximation technique} \label{alg5}
\begin{algorithmic}[1]
    
\Require
A tessellation $\mathcal{T}_h$ of $D$, stochastic diffusion coefficient $a(\bm{x},\omega)$, the initial control guess $\bm{f}_h^{(0)}$, tolerance $\epsilon$, and dimension reduction ratio $\tau$.
\Ensure
Numerical optimal control $\bm{f}_h^{*}$, and the expectation of optimal state $\bm{\mu}_{h,M,\tau}^*$.

\State Construct $V_h \in H_0^1(D)$ and sample the random field in diffusion coefficient $\widetilde{a}^m, m = 1,...,M$.

\State Assemble the stiffness matrices $\overline{\mathbb{A}}, \{\widetilde{\mathbb{A}}^m\}_{m=1}^M$ by Eq. (\ref{eq4.8}), and deal with the Dirichlet boundary conditions to $\overline{\mathbb{A}}$.

\State Determine the rank $k = \lceil \tau N \rceil$ and compute the sample states $\bm{u}_{h,\tau}^m$ using Algorithm \ref{alg3}.

\While{ $\Vert \nabla \widehat{\mathcal{J}}(\bm{f}_h^{(k)}) \Vert \: > \: \epsilon$ }

\State Compute the gradient $\nabla \widehat{\mathcal{J}}(\bm{f}_h^{(0)})$ or the hessian matrix $\Delta \widehat{\mathcal{J}}$ with respect to $\bm{f}_h^{(0)}$ according to Eq. (\ref{eq4.30}) and (\ref{eq4.31}).

\State Choose a suitable step size $\alpha^{(k)} \: > \: 0$.

\State Compute the step $\delta \bm{f}_h^{(k)}$ based on the derivative information $\nabla \widehat{\mathcal{J}}(\bm{f}_h^{(0)})$ or $\Delta \widehat{\mathcal{J}}$, and the step size $\alpha^{(k)}$.

\State Update the control by $\bm{f}_h^{(k+1)} \leftarrow \bm{f}_h^{(k)} + \delta \bm{f}_h^{(k)}$.

\State $k = k + 1$.

\EndWhile

\For{$m = 1,...,M$}    

\State Solve the state $\bm{u}_{h,\tau}^m = Z^m \bm{f}_h^{(k)}$ by Eq. (\ref{eq4.23}).

\EndFor

\State \Return the numerical optimal control $\bm{f}_h^* = \bm{f}_h^{(k)}$ and state $\bm{\mu}_{h,M,\tau}^* = \frac{1}{M} \sum_{m=1}^M \bm{u}_{h,\tau}^m$.

\end{algorithmic}
\end{breakablealgorithm}

In the optimization loop of Line 6 - 8 in Algorithm \ref{alg5}, the control $\bm{f}_h$ is updated by

\begin{equation} \label{eq4.32}
\bm{f}_h^{(k+1)} \; - \; \bm{f}_h^{(k)} \: := \: \delta \bm{f}_h^{(k)} \: = \: - \; \alpha^{(k)} d^{(k)},
\end{equation}

\noindent where the positive scalar $\alpha^{(k)}$ is the step length of the iteration, and $d^{(k)}$ denoted the descent direction. Note that the updating formula in Eq. (\ref{eq4.32}) allows various methods to obtain the descent direction and the step length. In the following numerical experiments, we apply multiple kinds of gradient-based unconstrained optimization methods, the trust-region method and the different rules of the line search techniques to obtain the control difference $\delta \bm{f}_h^{(k)}$. 

\subsubsection{Numerical Experiments on Elliptic SOCP}

In the following, we examine the numerical performance of Algorithm \ref{alg5} in solving the optimal control problem Eq. (\ref{eq4.18}) governed by the elliptic SPDE in Eq. (\ref{eq4.19}). The numerical example has the following settings. Let the domain $D = [0,1]^2$ triangularized as shown in Figure \ref{fig1}. We choose the regularization parameter $\beta = 1 \times 10^{-4}$, the tolerance parameter $\epsilon = 1 \times 10^{-3}$, the MC sample size $M = 200$, and the maximum amount of iterations in the line search method with the Wolfe condition $it_{max} = 50$ for the numerical implementation. The random diffusion coefficient $a(\bm{x},\omega)$ has the similar form as Eq. (\ref{eq4.17}) where $\sigma(\bm{x}, \omega)$ obeys the uniform distribution on the interval of $[-1,1]$ for any fixed $\bm{x} \in D$. The desired state $U(\bm{x})$ of the stochastic optimal control problem is given by $U(\bm{x}) = \sin(2 \pi x) \sin(2\pi x)$.

\begin{figure}[ht]
\centering
\includegraphics[width=1 \textwidth]{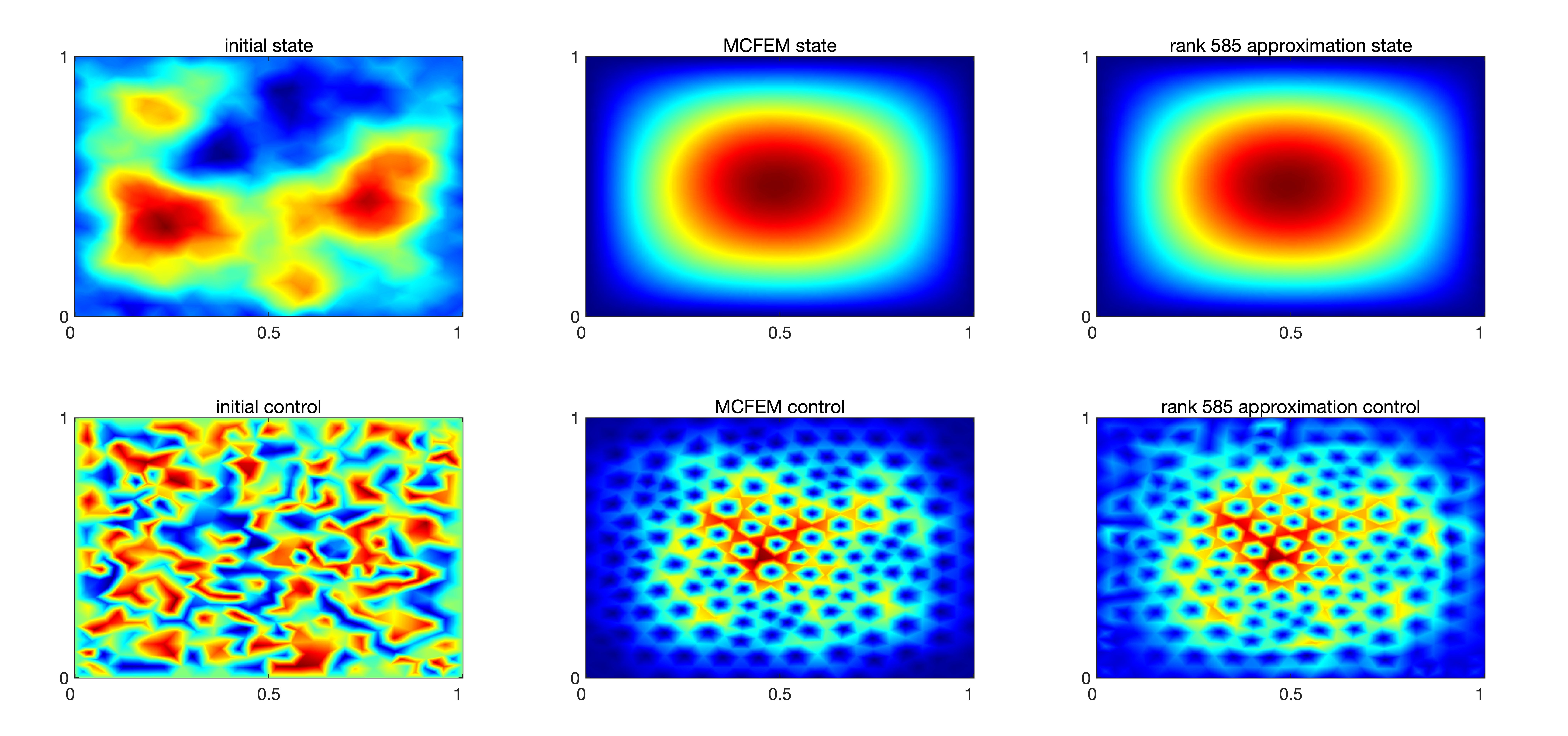}
\caption{The initial state $u^{(0)}(\bm{x})$ and control $f^{(0)}(\bm{x})$ (left), the reference state $u_{ref}(\bm{x})$ and control $f_{ref}(\bm{x})$ (middle), and the numerical optimal state $u^*(\bm{x})$ and control $f^*(\bm{x})$ (right).}
\label{fig4.5}
\end{figure}

We first examine the numerical performance of Algorithm \ref{alg5}. Here we set the dimension reduction rate $\tau = 0.88$ and apply the Newton's method with the line search method. The optimization loop is ended when the stopping criterion $\Vert \nabla \widehat{\mathcal{J}}(\bm{f}_h^{(k)}) \Vert \le 1 \times 10^{-3}$ is met. Meanwhile, the MCFEM approach ($\tau = 1$) is employed and its solutions are used as a reference. The numerical results of solving the stochastic elliptic optimal control problem in the above settings are demonstrated in the graphical form. Figure \ref{fig4.5} respectively presents the initial state and control, the reference solutions and the numerical optimal outputs. It is found from the plots that the state $u$ and control $f$ varies greatly after the optimization process, and their visualizations are similar to the reference solutions obtained by the MCFEM approach, which confirms the effectiveness of Algorithm \ref{alg5}.

\begin{figure}[ht]
\centering
\includegraphics[width=1 \textwidth]{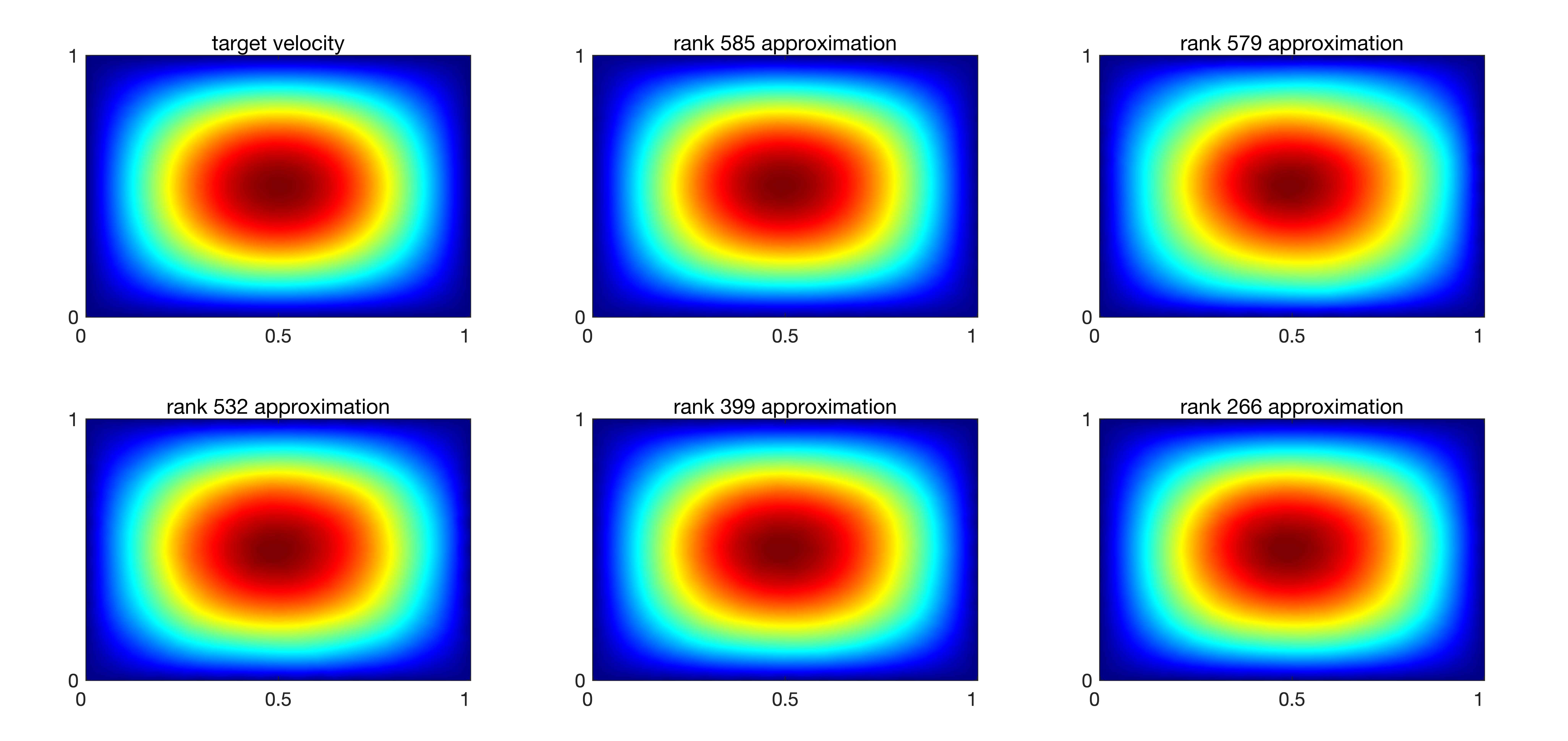}
\caption{The optimal state $u^*$ with 5 different dimension reduction ratio $\tau = 0.88$, $\tau = 0.87$, $\tau = 0.8$, $\tau = 0.6$, and $\tau = 0.4$.}
\label{fig4.6}
\end{figure}

Similarly we examine the sensitivity of Algorithm \ref{alg5} to the choice of the dimension reduction ratio $\tau$. Figure \ref{fig4.6} depicts the optimal state $u^*$ under 5 different
$\tau = 0.88, 0.87, 0.8, 0.6, 0.4$. We also present the comparison statistics in Table \ref{tab3}. The simulation outputs include the the CPU elapsed time, the errors between the optimal and desired state $\Vert u^{*} - U \Vert$, the final values of the gradient of the objective functional $\nabla\widehat{J}^{*}$, the initial and final Values of the objective functional, $\widehat{J}^{(0)}$ and $\widehat{J}^{*}$, and the ratios $\widehat{J}^{*}/\widehat{J}^{(0)}$. The columns of ratios and errors in Table \ref{tab3} measures the effect of the minimization procedure and the accuracy of numerical solutions. From the numerical results, the optimization accuracy of Algorithm \ref{alg5} increases with the grow of dimension reduction ratio $\tau$, while lower $\tau$ leads to less CPU elapsed time and storage. 

\begin{table}[ht]
\centering 
\caption{Simulation results for Different Dimension Reduction Ratios $\tau$ about the CPU elapsed time, Errors $\Vert u^{*} - U \Vert$, Final Value of Gradient of Objective Functional $\nabla\widehat{J}^{*}$, Initial and Final Value of Objective Functional, $\widehat{J}^{(0)}$ and $\widehat{J}^{*}$, and the Ratios $\widehat{J}^{*}/\widehat{J}^{(0)}$.}
\begin{tabular}{c|c|cccccccc}
\hline
\multicolumn{2}{c|}{} & Time & Error & $\nabla \widehat{J}^{*}$ & $\widehat{J}^{(0)}$ & $\widehat{J}^{*}$ & Ratio (\%) \\ \hline
\multicolumn{2}{c|}{Reference} & $0.9340$ & $1.1719$ & $0.000857506$ & $84.0994$ & $7.7592$ & $9.22$\\  \hline
\multirow{5}{*}{$\tau$} & $0.88$ & $0.7872$ & $1.2703$ & $0.000821937$ & $85.2419$ & $8.66532$ & $10.17$\\ 
& $0.87$ & $0.5613$ & $2.0302$ & $0.000658835$ & $84.5952$ & $13.7472$ & $16.25$\\ 
& $0.8$ & $0.5535$ & $2.0675$ & $0.000636466$ & $84.4904$ & $14.2412$ & $16.86$\\ 
& $0.6$ & $0.5437$ & $2.4200$ & $0.000576569$ & $84.5121$ & $20.2741$ & $19.63$\\ 
& $0.4$ & $0.5347$ & $2.5441$ & $0.000560877$ & $84.5693$ & $17.414$ & $20.59$\\ \hline
\end{tabular}
\label{tab3}
\end{table}

One of the advantages of Algorithm \ref{alg5} is that it can make use of various types of gradient-based algorithms to solve the unconstrained minimization problem in Eq. (\ref{eq4.26}). In the experiment, we employ the following five optimization algorithms: the steepest descent method (SDM), the stochastic gradient descent (SGD) method, the Newton's method, the BFGS method and the trust region with dogleg method (TRM). Table \ref{tab4} presents the comparison results, which shows that Algorithm \ref{alg5} has extensive suitability to unconstrained optimization methods. Moreover, the SGM does the best job in solving the stochastic optimal control problem among the five unconstrained optimization algorithms, which is not surprise, since SGM does not need to compute the gradients of all the MC realizations and significantly reduces the work complexity. As a result, the stochastic gradient method in a unique position when dealing with high-dimensional random space.

\begin{table}[ht]
\centering 
\caption{Simulation results for Different Unconstrained Optimization Methods about Iterations, the CPU elapsed time, Errors $\Vert u^{*} - U \Vert$, and the Ratios $\widehat{J}^{*}/\widehat{J}^{(0)}$.}
\begin{tabular}{cccccc}
\hline
Method & SDM & SGD & Newton & BFGS & TRM \\ \hline
Iteration & $75$ & $94$ & $8$ & $28$ & $41$ \\
Time & $5.8367$ & $0.4774$ & $0.7872$ & $1.5499$ & $0.7165$ \\
Error & $3.8084$ & $0.9431$ & $1.2703$ & $1.7866$ & $0.9221$ \\
Ratio $\widehat{J}^{*}/\widehat{J}^{(0)}$ (\%) & $16.19$ & $7.12$ & $10.17$ & $12.15$ & $9.82$\\ \hline
\end{tabular}
\label{tab4}
\end{table}

\section{Conclusions and Discussions}

In this article, we propose a fast algorithm for efficiently solving partial differential equations with perturbations, which have probabilistic discretized formulations of $\mathbb{A}_m \bm{u}_m = \bm{b}$. By splitting the perturbation and applying a novel low rank approximation method for the collection of large-scale perturbed matrices $\{\mathbb{A}_m\}_{m=1}^M$, our algorithm can significantly reduce the complexity and storage of the computation of matrix inversion. To demonstrate the versatility and applicability of our algorithm, we apply it to address two crucial applications: stochastic elliptic partial differential equations and stochastic optimal control problems governed by elliptic PDE constraints. Depending on different dimension reduction ratio $\tau$, our algorithm can construct a high precision numerical solution or roughly depict the main sketch of the exact solution as a pre-processing. Meanwhile, the algorithm makes good use of various types of gradient-based unconstrained optimization methods, which makes it more attractive for large-scale problems from a computational point-of-view. Numerical results from both applications validate the feasibility and the effectiveness of the proposed algorithm.

However, there are still many open questions to be answered. Firstly, future research should consider the negative effects from the instability of perturbed linear systems, and some pre-processing measures for ill-conditioned stiffness matrices are necessary. The algorithm presents good effectiveness in solving linear elliptic SPDE and SOCP. However, the computational complexity and requirement levels up for unsteady or nonlinear systems, such as the Navier-Stokes equation. Therefore, the efficiency and feasibility of our algorithm needs to be further studied. In addition, we can make a further investigation on the trade-off between the data compression ratio and the numerical precision in the low rank matrix approximation process. To reduce the matrix reconstruction errors, future works may aim at the higher-degree formulations, such as $\min \Vert \mathbb{A} - \mathbb{U} \mathbb{U}^T \mathbb{V} \mathbb{V}^T\Vert$. On the other hand, pre-processing of the perturbed matrices $\widetilde{\mathbb{A}}_m$ might be a potential research topic. We could pre-classify the matrices in a clustering manner and use different $\mathbb{U}, \mathbb{V}$ in the matrix approximation. We hope that in the near future, we will be able to find answers to these questions.

\section*{Acknowledgments}
The authors would like to thank the anonymous referees and the editor for their valuable comments and suggestions, which led to considerable improvement of the article.

\section*{Conflict of Interest}
All authors declare that there are no conflicts of interest regarding the publication of this paper.

\bibliographystyle{elsarticle-harv}
\bibliography{bib}

\end{document}